\documentclass[a4paper,12pt] {article}
\usepackage{graphics}
\usepackage{graphicx}
\usepackage[english]{babel}
\usepackage{array}
\usepackage{multirow}
\usepackage{mathrsfs}
\usepackage{amsmath}
\usepackage{amssymb}
\usepackage{float}
\usepackage[position=top,singlelinecheck=false,justification=raggedleft,captionskip=0pt,topadjust=-8pt]{subfig}
\usepackage{multirow}
\usepackage{xcolor}
\usepackage{placeins}
\usepackage{cite}
\usepackage[abs]{overpic}
\usepackage{caption}
\usepackage{stix}

\definecolor{gray}{RGB}{153,153,153}
\definecolor{blue}{RGB}{0,114,189}
\definecolor{orange}{RGB}{222,125,0}
\definecolor{green}{RGB}{153,204,0}
\definecolor{purple}{RGB}{126,47,142}
\definecolor{pink}{RGB}{218,179,255}
\definecolor{lblue}{RGB}{173,235,255}
\definecolor{ocean}{RGB}{0,191,191}
\definecolor{myred}{RGB}{219,75,81}
\definecolor{red}{RGB}{255,0,0}
\definecolor{brown}{RGB}{153,51,0}
\definecolor{dark_green}{RGB}{119,172,48}

\begin{document}
\title{\Large\textbf{Force appropriation of nonlinear structures}}
\author{\textbf{L.~Renson, T.L.~Hill, D.A.~Ehrhardt, D.A.W.~Barton, S.A.~Neild} \vspace{2mm}\\
\normalsize{Faculty of Engineering, University of Bristol, UK.}\\
\normalsize{Corresponding author: l.renson@bristol.ac.uk}}

\date{\normalsize{submitted \textit{Dec. 2017}}}
\maketitle

\begin{abstract}
Nonlinear normal modes (NNMs) are widely used as a tool for developing mathematical models of nonlinear structures and understanding their dynamics. NNMs can be identified experimentally through a phase quadrature condition between the system response and the applied excitation. This paper demonstrates that this commonly-used quadrature condition can give results that are significantly different from the true NNM, in particular when the excitation applied to the system is limited to one input force, as is frequently used in practice. The system studied is a clamped-clamped cross beam with two closely-spaced modes. This paper shows that the regions where the quadrature condition is (in)accurate can be qualitatively captured by analysing transfer of energy between the modes of the system, leading to a discussion of the appropriate number of input forces and their locations across the structure.\\

\end{abstract}

\section{Introduction}
The dynamics of a nonlinear structure is often characterised by the presence of interactions between the different degrees of freedom (DOFs) of the system. Such interactions arise from the presence of nonlinear coupling terms in the equations of motion, which, contrary to the case of linear systems, cannot be completely removed by a transformation of coordinates that results in a set of independent oscillators. In structural dynamics, DOFs can be associated to particular modes of vibration. Interactions between modes can lead to energy transfer between different parts of a structure, which can jeopardise structural integrity~\cite{Noel14}.
Conversely, modal interactions can be exploited to improve design performance, for example, in vibration absorbers~\cite{VakakisTET}, high-performance frequency dividers~\cite{Qalandar14}, as well as a number of micro- and nano-mechanical devices~\cite{Rhoads08}. Investigating the presence of such interactions is therefore key for understanding the dynamic behaviour of many nonlinear structures.

The theory of nonlinear normal modes (NNMs) has been successfully used to analyse and predict the presence of modal interactions in nonlinear structures\cite{Renson15,Ducceschi13,Liu16}. NNMs can be defined as families of periodic oscillations of a conservative (i.e. undamped) system~\cite{Kerschen09}. NNMs generally do not have the convenient mathematical properties that linear modes have, such as orthogonality~\cite{GeradinBook}; nevertheless, they still correspond to invariant properties which can be used to interpret a number of important nonlinear phenomena, such as localisation~\cite{Vakakis97}, mode bifurcations~\cite{Sarrouy11}, energy transfer~\cite{Georgiades09b,Vaurigaud11}, isolated responses~\cite{Kuether15,Hill16b} and modal interactions. NNMs can also trace the evolution of the resonance frequency of a damped, harmonically-forced system~\cite{Kerschen09}, which is valuable from an engineering perspective because this is where displacements are often maximum and the structure is at the greatest risk of failure. 

The relationship between NNMs, which are conservative properties, and the response of a damped, forced system occurs when the phases of the external input forces and the phases of the responses of the structure are in quadrature, i.e. when they have a phase difference of $90$ degrees. At quadrature, external input forces exactly counterbalance the internal damping forces and hence the system responds as the underlying conservative system~\cite{Peeters11b}. This phase quadrature condition is very attractive because it enables the identification of the NNMs of a system through the experimental testing of the physical structure~\cite{Peeters11c,Renson16,Peter17}. Identified NNMs can, in turn, be exploited for parameter estimation~\cite{Hill16} and damage detection \cite{Lacarbonara16}, compared to theoretical predictions for model updating and validation~\cite{Peter15,Song18,Ehrhardt17} or, simply, for monitoring the evolution of the resonance frequency and amplitude of the response.

The mathematical derivation of the phase quadrature condition usually relies on three major assumptions. First, the damping in the structure is considered to be linear (although some forms of nonlinear damping can be considered~\cite{Renson16}). Second, the applied excitation is multi-harmonic and each of its harmonics is in phase quadrature with the corresponding harmonic in the response. Finally, the excitation is applied to all the DOFs of the structure. For structures that behave essentially like SDOF oscillators with small modal couplings, a single-point, single-harmonic excitation may provide an accurate estimation of the NNMs. However, this is no longer the case for structures where strong couplings between modes exist. For instance, in Ref.~\cite{Ehrhardt16}, a third harmonic component was added to the applied excitation in order to capture the NNMs of a doubly-clamped beam featuring a 3:1 modal interaction.

\begin{figure}[t]
\begin{center}
\includegraphics[width=0.95\textwidth]{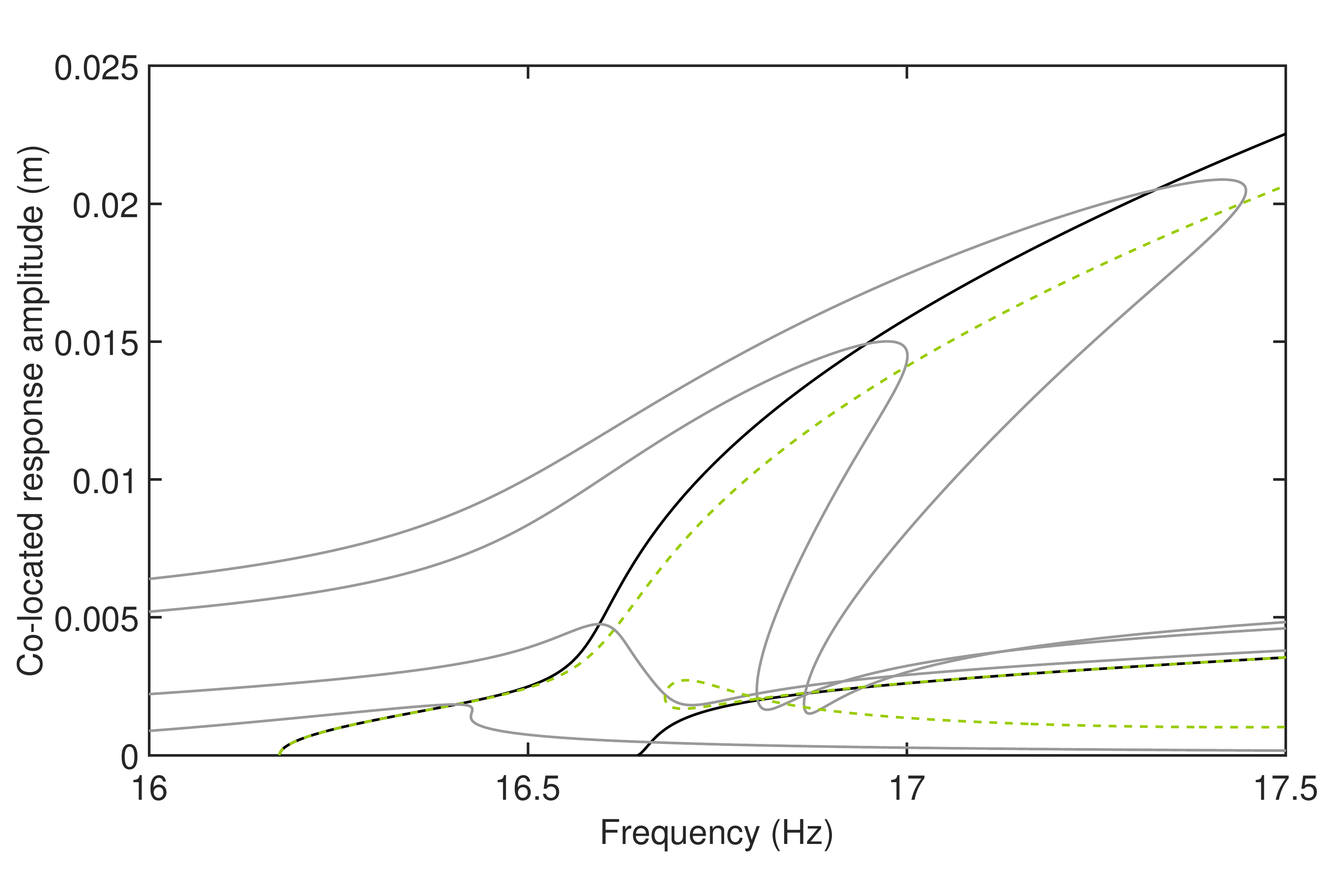}
\caption{Nonlinear frequency response of a clamped-clamped cross beam structure for four different levels of excitation (\textcolor{gray}{$\boldsymbol{-}$}). The response of the structure is shown at the excited DOF (co-located). The green curves (\textcolor{green}{$\boldsymbol{--}$}) trace out the loci of periodic responses that are in phase quadrature with the excitation, and the thick black curves (\textcolor{black}{$\boldsymbol{-}$}) show the NNMs of the underlying conservative system.}
\label{fig:11beam_motiv}
\end{center}
\end{figure}

In this paper, we investigate the case of a structure with two closely-spaced modes, leading to a 1:1 modal interaction. This is representative of a number of structural applications, including complex industrial structures with high modal density~\cite{Kerschen13} and structures with (broken) symmetries~\cite{Touze02}. It is shown that the spatial distribution of the excitation has a strong impact on phase quadrature, which may no longer correspond to the NNMs of the underlying conservative system. Figure~\ref{fig:11beam_motiv} illustrates this issue. The frequency response of the example structure (\textcolor{gray}{$\boldsymbol{-}$}) is shown for four different amplitudes of single-point, single-harmonic excitation. The NNMs of this system are represented by solid lines (\textcolor{black}{$\boldsymbol{-}$}). It can be seen that the resonant peaks of two high-amplitude frequency responses do not correspond to the first NNM. Furthermore, the point of phase quadrature (\textcolor{green}{$\boldsymbol{--}$}), which varies both in frequency and amplitude, diverges from this NNM. This error appears to be difficult to predict based on the low-amplitude behaviour as the agreement between the quadrature and the first NNM is almost perfect at low amplitude. Phase quadrature points for the second NNM exist only for sufficiently-large response amplitudes and form an alpha-shaped curve, which only partly corresponds to the NNM. This demonstrates that, in this case, the quadrature curves are not analogous to the NNMs of the system. Even larger discrepancies will be observed and discussed in this paper.

The experimental identification of NNMs is usually carried out using a single exciter. The objective of this paper is to investigate the discrepancies that can arise between quadrature curves and NNMs in such excitation conditions. The example structure will be exploited to demonstrate the potential consequences of inappropriate forcing and study the influence of the spatial position of the excitation. Based on energy arguments, this paper will also demonstrate that inappropriate excitations lead to energy transfer and phase differences between DOFs, which are responsible for inaccurate quadrature conditions. In Section~\ref{sec:11beam}, the example structure and its NNMs are discussed. Then in Section~\ref{sec:theory}, the excitation necessary to isolate a specific NNM motion is calculated using energy arguments and is numerically demonstrated for the theoretical case in which an arbitrary number of DOFs can be excited. The case where the external excitation is restricted to a single input force is then discussed in Section~\ref{sec:nnm_1fphys}. Section~\ref{sec:pred_energy} further develops the energy arguments of Sections~\ref{sec:theory} to analyse energy transfer between modes and qualitatively identify the regions where quadrature is (in)accurate. The approach is exploited to discuss the merits of different excitation locations. Conclusions of this study are drawn in Section~\ref{sec:conclu}.

\section{Example structure---a beam with two closely-spaced modes}\label{sec:11beam}
\begin{figure}[th]
\centering
\includegraphics[width=0.5\textwidth]{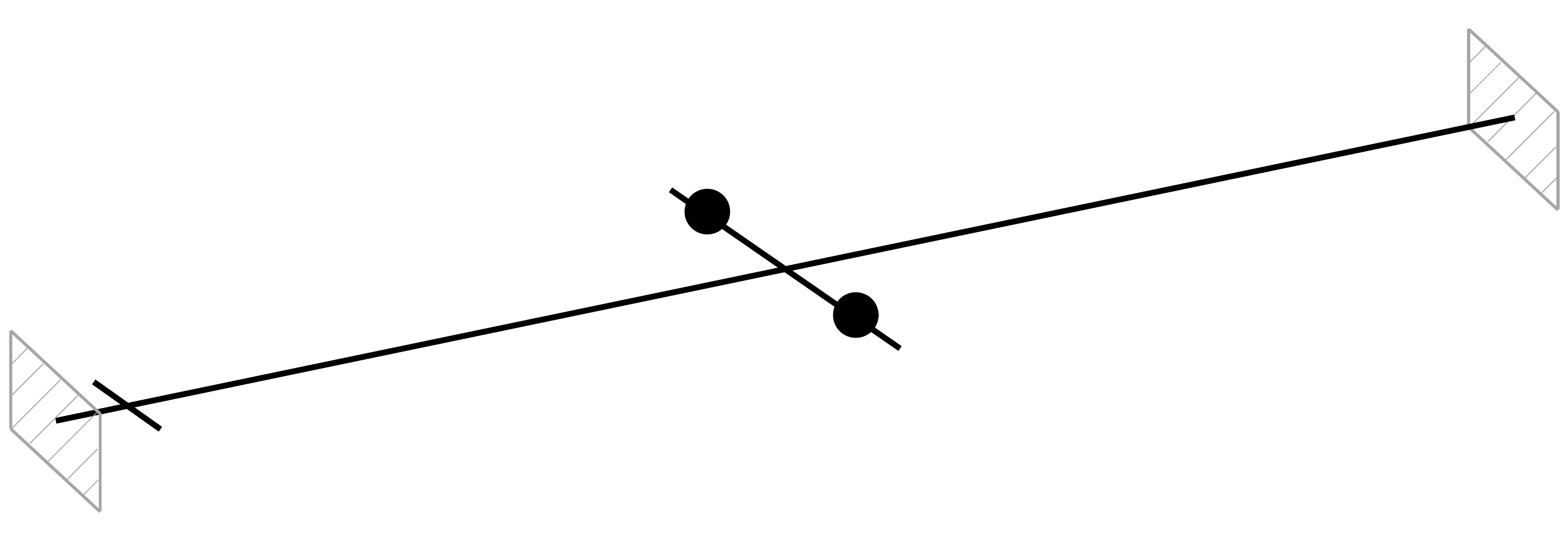}
\caption{Schematic of the clamped-clamped cross beam structure. The two masses are asymmetrically positioned on the cross beam, and tuned such that the first two linear natural frequencies are $0.5$\,Hz apart. A smaller cross beam, without masses, provides additional locations for external excitation.}
\label{fig:rig}
\end{figure}

The cross beam structure, illustrated in Figure~\ref{fig:rig} and studied previously in \cite{Ehrhardt17sub},
is used as a motivating example throughout this paper. The doubly-clamped boundary conditions couple the transverse (bending) and axial (stretching) motions of the main beam, leading to nonlinear geometric effects~\cite{WaggNeildBook}. Such nonlinearity is reminiscent of a number of applications in aeronautics and micro-mechanics where nonlinearity arises from finite displacements~\cite{Silva94,Westra12,Polunin17}. Two masses are attached to the cross beam in order to provide a means of adjusting both the torsional inertia and the symmetry of the system. When the system is perfectly symmetric, the first two linear modes correspond to pure bending and pure torsion motions of the main beam. However, in this paper, the masses are asymmetrically positioned --- one mass is approximatively 9 mm closer to the mean beam compared to the other one --- and adjusted such that the first two natural frequencies are close. This leads to bending and torsion components in both modes, although the first mode is dominated by bending and the second by torsion as the break in symmetry is only slight.

The structure is modelled as a two-degree-of-freedom, nonlinear, modal model capturing the dynamics of the first two modes of interest. Here we use the term mode with reference to a linearised version of the model, hence the modal model may, and does, exhibit mode cross-coupling due to nonlinearities. The equations of motion take the form 
\begin{equation}
\mathbf{\ddot{q}}(t) + \mathbf{\Xi} \mathbf{\dot{q}}(t)  + \boldsymbol{\Lambda} \mathbf{q}(t) + \mathbf{N}_{\mathbf{q}}(\mathbf{q}(t)) = \mathbf{p}(t),
\label{eq:modal_damped}
\end{equation}
where $\boldsymbol{\Lambda}$ is the linear modal stiffness matrix, $\mathbf{\Xi}$  is a diagonal matrix, $\mathbf{\Xi}_{i,i} = 2 \zeta_i \omega_{ni}$, $i = 1, 2$, representing linear modal damping and $\mathbf{p}(t)$ is the modal forcing vector. The vector $\mathbf{N}_{\mathbf{q}}(\mathbf{q})$ of nonlinear internal forces is chosen to be a sum of quadratic and cubic terms as
\begin{equation}
\mathbf{N}_{\mathbf{q}}(\mathbf{q}) = \left( \begin{array}{c}
	\alpha_1 q_1^2 + 2 \alpha_2 q_1 q_2 + \alpha_3 q_2^2 + \gamma_1 q_1^3 + 3\gamma_2 q_1^2 q_2 + \gamma_3 q_1 q_2^2 + \gamma_4 q_2^3 \\
	\alpha_2 q_1^2 + 2 \alpha_3 q_1 q_2 + \alpha_4 q_2^2 + \gamma_2 q_1^3 + \gamma_3 q_1^2 q_2 + 3 \gamma_4 q_1 q_2^2 + \gamma_5 q_2^3 
\end{array}\right).
\label{eq:nq_rom}
\end{equation}

This modal model was obtained through the reduction of a large-scale finite-element (FE) model of the structure using the implicit condensation and expansion method (ICE)~\cite{Kuether15b}. The FE model was built in Abaqus using nonlinear geometric elements and comprised of 1914 DOFs. The equations of motion of this physical model take the form
\begin{equation}
\mathbf{M} \ddot{\mathbf{x}}(t) + \mathbf{K} \mathbf{x}(t) + \mathbf{f}_{\text{nl}}(\mathbf{x}(t)) = \mathbf{f}(t),
\label{eq:eom}
\end{equation}
where $\mathbf{M}$ and $\mathbf{K}$ are the linear mass and stiffness matrices, respectively, $\mathbf{f}_{\text{nl}}(\mathbf{x}(t))$ is the nonlinear force vector and $\mathbf{f}(t)$ is the vector of external forces applied to the structure. The modal forcing vector $\mathbf{p}(t)$ is obtained by projection of the physical forces $\mathbf{f}(t)$ onto the reduced basis of mass-normalised linear mode shapes $\mathbf{\Phi}$. Note that only external excitation forces applied to the vertical DOFs of the structure will be considered in this paper. The transformation from the reduced modal space to the full physical space is given by $\mathbf{x} = \mathbf{\Phi} \mathbf{q}$. The linear and nonlinear properties of the system are listed in Table~\ref{tab:11beam_all_prop}. Damping parameters were taken from Refs.~\cite{Ehrhardt17sub,RensonISMA2016} where a similar system was investigated experimentally. 
\begin{table}
\centering
\begin{tabular}{cc|cc|cccc|ccccc}
$\omega_{n1}$ & $\omega_{n2}$ & $\zeta_{1}$ & $\zeta_{2}$ & $\alpha_1$ & $\alpha_2$ & $\alpha_3$ & $\alpha_4$ & $\gamma_1$ & $\gamma_2$ & $\gamma_3$ & $\gamma_4$ & $\gamma_5$ \\
\hline
\multicolumn{2}{c}{$(\mbox{rad} \; \mbox{s}^{-1})$} & \multicolumn{2}{|c|}{$(\times 10^{-3})$} & \multicolumn{4}{c|}{$(\times 1)$} & \multicolumn{5}{c}{$(\times 10^6)$}\\
101.61 & 104.58 & 7.6 & 2.6 & 56.7 & -52.4 & -14.9 & 42.7 & 128 & 32 & 25 & 2 & 0.8 
\end{tabular}
\caption{Linear and nonlinear model parameters of the clamped-clamped cross beam structure.}
\label{tab:11beam_all_prop}
\end{table}

The first two NNMs, or modes of the nonlinear system, for the beam structure were calculated using the harmonic balance technique. Assuming that the responses of the modes may be written as
\begin{equation}
	q_i \approx u_i = U_i \cos\left(\Omega t - \phi_i \right),
\label{eq:assumed_sol}
\end{equation}
where $U_i$ and $\phi_i$ denote the fundamental response amplitude and phase of the $i^{\mbox{\scriptsize th}}$ mode, respectively and $\Omega$ represents the response frequency, the equations governing the frequency-amplitude dependence of the first and second NNMs may be written as 
\begin{subequations}
\begin{eqnarray}
	\left(
		\omega_{n1}^2 - \Omega^2
	\right) U_1
	+ \dfrac{3 \gamma_1}{4} U_1^3
	+ p \dfrac{3 U_2}{4} \left[
		3 \gamma_2 U_1^2
		+ \gamma_4 U_2^2
	\right]
	+ \dfrac{3 \gamma_3}{4} U_1 U_2^2
	&\hspace{-6pt}=\hspace{-6pt}& 0,
\label{eq:resonant_sol1}%
\\
	\left(
		\omega_{n2}^2 - \Omega^2
	\right) U_2
	+ p \dfrac{3 U_1}{4} \left[
		\gamma_2 U_1^2
		+ 3 \gamma_4 U_2^2
	\right] 
	+ \dfrac{3 \gamma_3}{4} U_1^2 U_2
	+ \dfrac{3 \gamma_5}{4} U_2^3
	&\hspace{-6pt}=\hspace{-6pt}& 0,
\label{eq:resonant_sol2}%
\end{eqnarray}%
\label{eq:resonant_sols}%
\end{subequations}%
where ${ p = -1 }$ when the modes are in anti-phase (first NNM) and ${ p = +1 }$ when the modes are in-phase (second NNM). Details of this derivation can be found in Appendix~\ref{app:NNMs}. The NNMs were also obtained numerically using an algorithm combining a harmonic balance solution and pseudo-arclength continuation, but other methods could have been used~\cite{RensonReview}. Within the range of frequencies considered, the agreement between analytical and numerical results is excellent. 

\begin{figure}[tbp]
\centering
\captionsetup[subfloat]{justification=raggedright}
\subfloat[]{\includegraphics[width=1\textwidth]{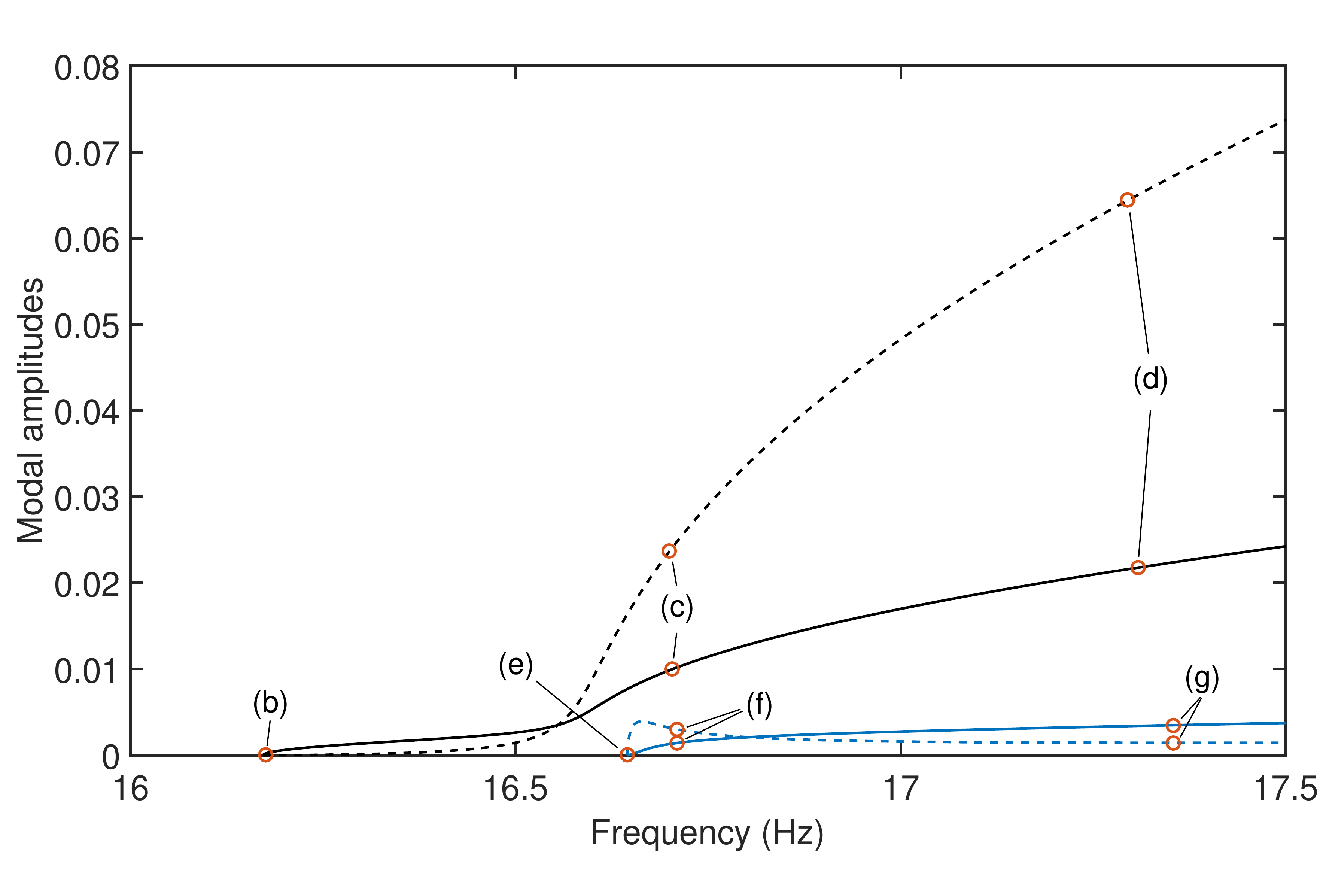}}\\
\captionsetup[subfloat]{justification=raggedright,topadjust=0pt}
\begin{tabular*}{0.9\textwidth}{@{\extracolsep{\fill}} c c c}
\subfloat[]{\label{shape1}\includegraphics[width=0.3\textwidth]{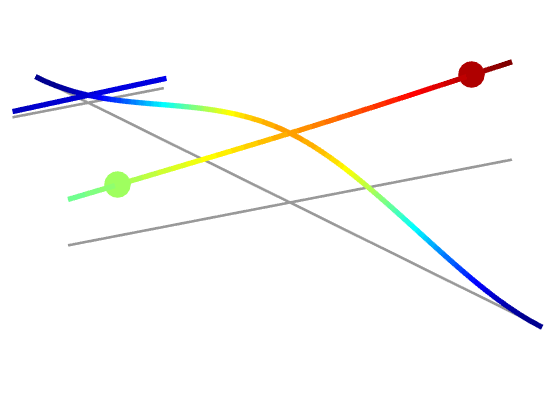}} &
\subfloat[]{\label{shape2}\includegraphics[width=0.3\textwidth]{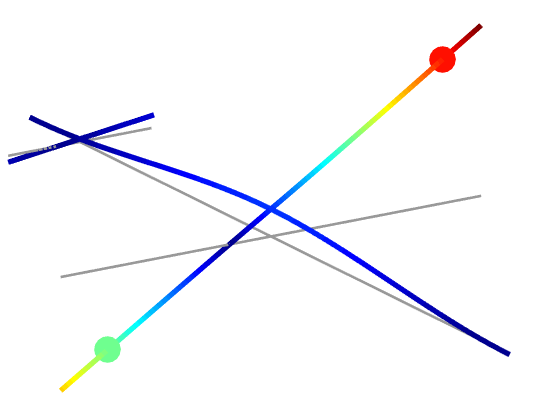}} &
\subfloat[]{\label{shape3}\includegraphics[width=0.3\textwidth]{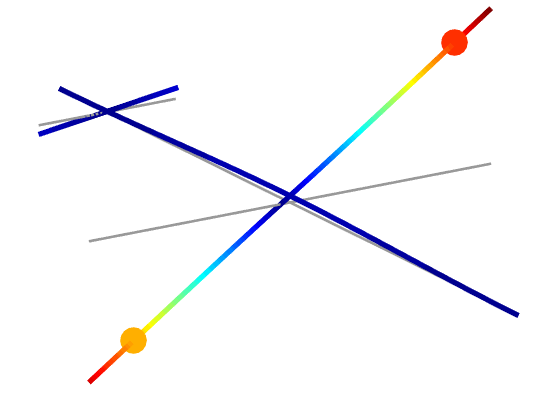}}\\
\subfloat[]{\label{shape4}\includegraphics[width=0.3\textwidth]{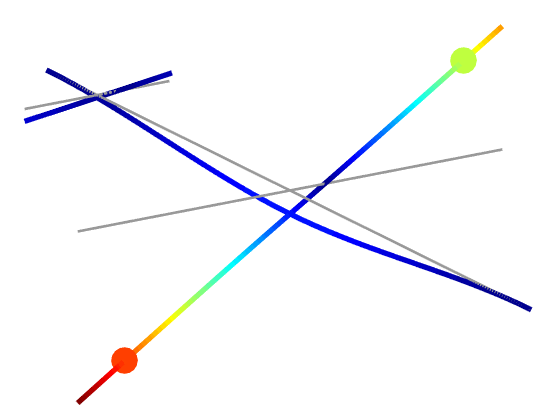}} &
\subfloat[]{\label{shape5}\includegraphics[width=0.3\textwidth]{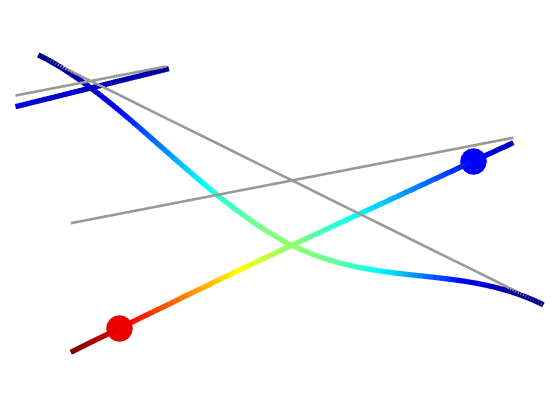}} &
\subfloat[]{\label{shape6}\includegraphics[width=0.3\textwidth]{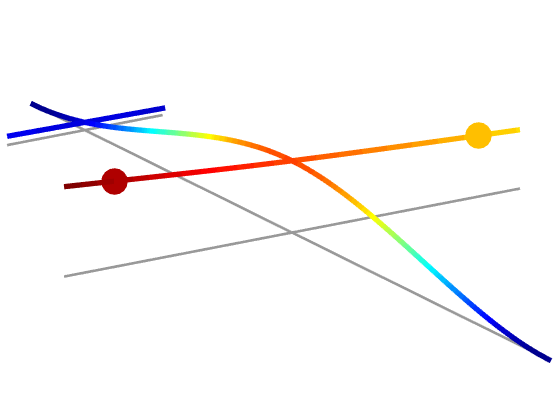}}\\
\end{tabular*}
\caption{NNMs of the beam structure. NNM 1: amplitude of modes 1 (\textcolor{black}{$\boldsymbol{-}$}) and 2 (\textcolor{black}{$\boldsymbol{--}$}). NNM 2: amplitude of modes 1 (\textcolor{blue}{$\boldsymbol{-}$}) and 2 (\textcolor{blue}{$\boldsymbol{--}$}). (b-g) Deformation of the structure  at the locations reported on the frequency-amplitude curves.}
\label{fig:nnms}
\end{figure}

The NNMs of the cross beam structure are presented in Figure~\ref{fig:nnms}, giving the amplitude of the first (solid) and second (dashed) linear modes as a function of the oscillation frequency. Due to the proximity of the linear modes and the break in symmetry, the frequency of the NNMs approach each other before veering off. This veering phenomenon is commonly observed in near-symmetric structures with closely-spaced modes~\cite{Lacarbonara05}. The shapes of the response of the structure vary significantly with the response amplitude and lead to a swap of mode shapes between NNMs from low (Figures~\ref{fig:nnms}(b, e)) to high amplitude (Figures~\ref{fig:nnms}(d, g)). Note the change of phase between the first and second NNM (compare, for instance, Figures~\ref{fig:nnms}(b) and~\ref{fig:nnms}(g)).

\section{Phase quadrature---arbitrary excitation conditions}\label{sec:theory}
The extension of phase quadrature to nonlinear systems was theoretically investigated in Ref.~\cite{Peeters11b}. Following the principles of linear force appropriation~\cite{Wright99}, several methods~\cite{Peeters11c,Renson16,Peter17} have exploited phase quadrature to identify NNMs experimentally. The structures investigated so far have, however, mostly been limited to systems that have, or can be approximated to, one DOF~\cite{Zapico13,Peeterse11c}. The objective of this section is to numerically demonstrate that the phase quadrature condition can still be used in the presence of strong nonlinear modal couplings, provided that a sufficient number of external forces are applied to the systems. However, such excitations may not be considered in practice, even for multi-DOF structures, as they are more difficult to implement experimentally. As such, the case where an insufficient number of forces are used will be considered in Section~\ref{sec:nnm_1fphys}. This section also introduces the energy balance approach and exploits it to determine the excitation amplitude required to obtain a specific NNM motion. This approach will be revisited in Section~\ref{sec:pred_energy} to analyse the accuracy of quadrature curves.

\subsection{Two-point excitation}\label{sec:11Exact}
\begin{figure}[th]
\centering
\begin{tabular*}{0.95\textwidth}{@{\extracolsep{\fill}} c c}
\subfloat[]{\label{flnm1}\includegraphics[width=0.45\textwidth]{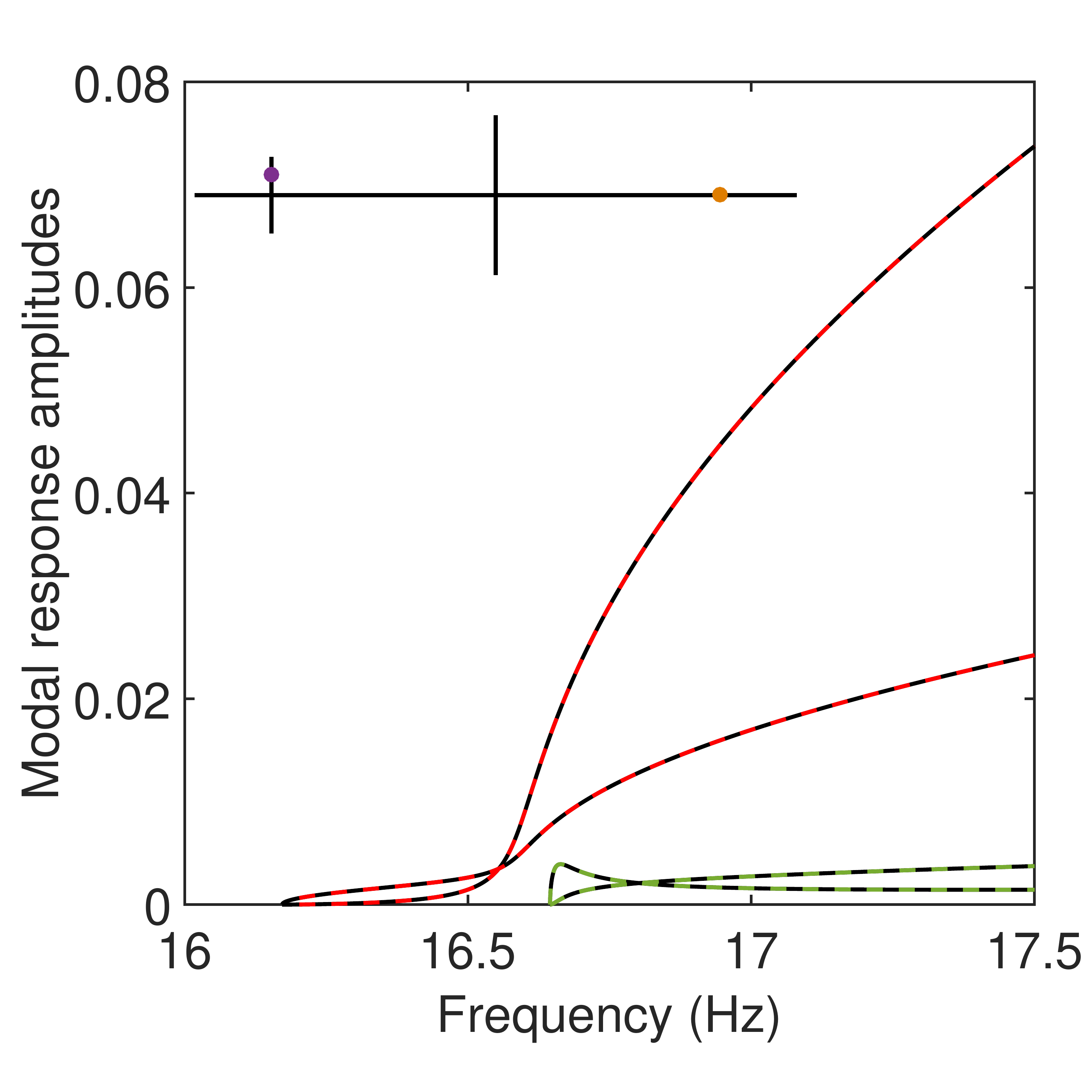}}&
\subfloat[]{\label{flnm2}\includegraphics[width=0.45\textwidth]{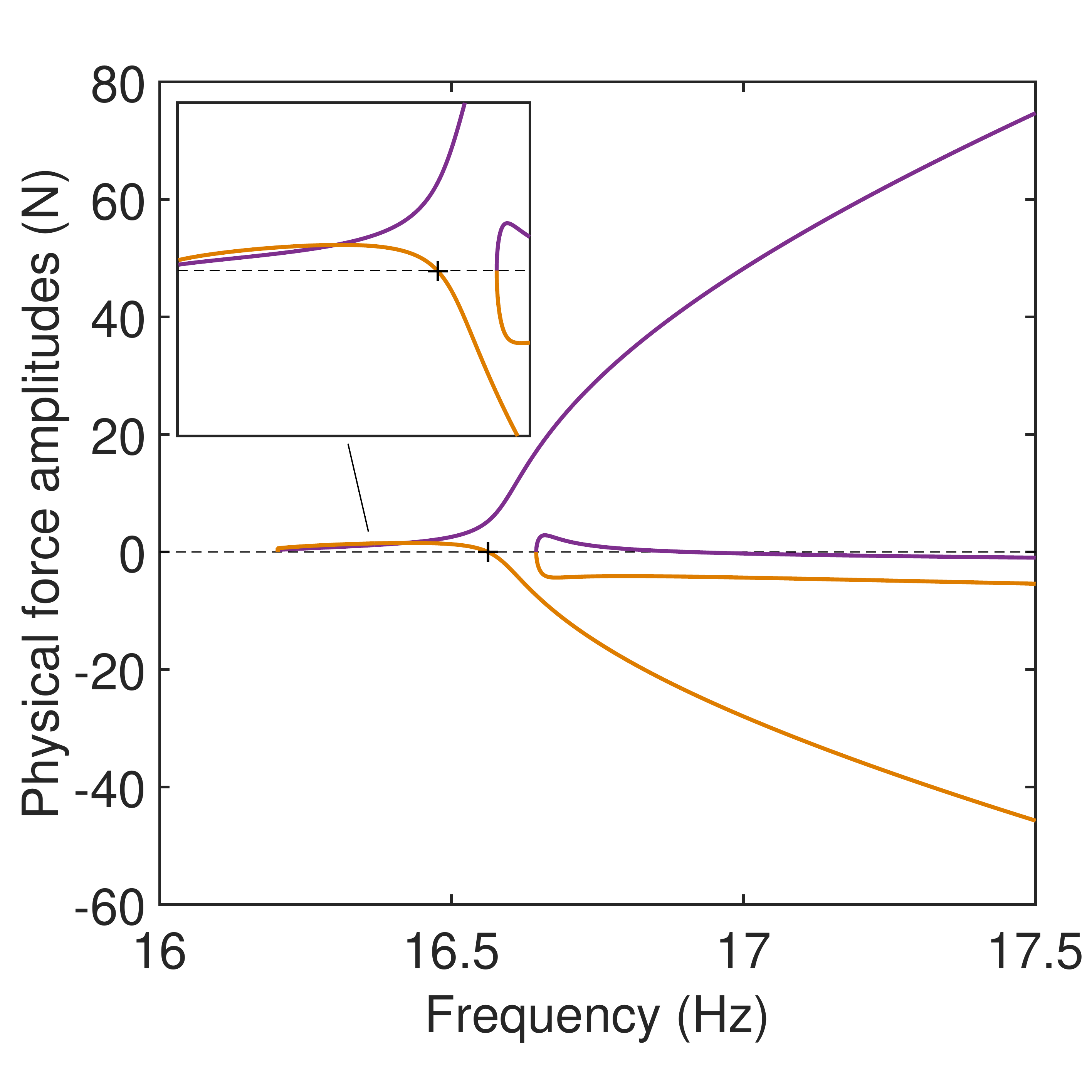}}
\end{tabular*}
\caption{(a) Comparison between the NNMs (\textcolor{black}{$\boldsymbol{-}$}) and the phase quadrature curves (\textcolor{red}{$\boldsymbol{--}$},\textcolor{dark_green}{$\boldsymbol{--}$}) obtained by applying a single-harmonic excitation at the locations marked in the structure schematic (\textcolor{purple}{$\boldsymbol{\bullet}$},\textcolor{orange}{$\boldsymbol{\bullet}$}). (b) Amplitudes of the physical forces required to isolate the first and second NNMs. Force amplitude curves are coloured according to their respective excitation locations.}
\label{fig:flnm12_nnms_freeratio_phys}
\end{figure}

The beam structure is excited with two synchronous, sinusoidal forces --- one on the small cross beam and one on the main beam (\textcolor{purple}{$\boldsymbol{\bullet}$} and \textcolor{orange}{$\boldsymbol{\bullet}$} in Figure~\ref{fig:flnm12_nnms_freeratio_phys}(a)). Starting from low-amplitude responses, the periodic responses that are in quadrature with the excitation are found and their evolution for increasing forcing amplitudes traced using numerical continuation. Two quadrature conditions are defined using the phase between the excitation force and the response at the excitation point (co-located). Note that the reduced-order modal model, Eq.~\eqref{eq:modal_damped}, is used for this simulation, and hence a modal transformation is required to extract the physical displacements.

The resulting quadrature curves are shown in Figure~\ref{fig:flnm12_nnms_freeratio_phys}(a). A near-perfect agreement with the NNMs is observed because the forcing required to isolate the NNMs can be reproduced with two input forces. This may appear obvious as Eq.~\eqref{eq:modal_damped} is reduced to two modes. However, this result still holds for a modal model comprising of five modes. Whilst some coupling between the first two modes and the three additional (bending) modes exists, it is negligible --- the amplitudes of the three higher modes are at least three orders of magnitude smaller than first two modes. As such, the NNMs are still very well-approximated by the quadrature curves obtained with two input forces. Note that the higher harmonics are not considered in this study because they are negligible compared to the fundamental component.

Figure~\ref{fig:flnm12_nnms_freeratio_phys}(b) shows the excitation amplitudes at the small cross and main beam required to isolate both NNMs. For the first NNM, both excitations have similar amplitudes although the excitation on the main beam is slightly larger. At approximately 16.55 Hz, the amplitude of the excitation applied to the main beam changes sign to become negative. When the excitation amplitude at the main beam crosses zero ($\boldsymbol{+}$), the excitation applied to the small cross alone provides the excitation necessary to perfectly isolate the first NNM. This observation will be further discussed in Section~\ref{sec:nnm_1fphys}.

As the excitation amplitude changes sign, it also changes its phase from $0$ to $-\pi$ rad. However, the phase of the response at the excitation location remains unchanged and the excitation now exhibits a $\pi/2$-phase lag with the response. This may appear surprising because, for single-DOF systems, the excitation is necessarily leading the response (provided damping is positive). In fact, past 16.55 Hz, the excitation applied to the small cross alone (\textcolor{purple}{$\boldsymbol{\bullet}$}) leads to an excessive excitation of the modal coordinates. The excitation applied to the main beam (\textcolor{orange}{$\boldsymbol{\bullet}$}), being negative, moderates this, and so achieving the forcing required to isolate the NNMs. Note that the modal forces resulting from the sum of the physical excitations still shows a $\pi/2$-phase lead with the response, which shows that the overall role of the excitation forces applied to the system remains to compensate internal damping. This observation will now be examined further using an energy-based approach.

\subsection{Energy balance}\label{sec:balance}
Given a specific NNM motion, it is possible to calculate the excitation forces required to isolate its motion using, for instance, harmonic balance~\cite{Peeters11b}. An alternative approach, which is considered here, is to use energy arguments as introduced in~\cite{Hill2014b} and later developed in \cite{Hill15}. This approach is based on the observation that, for any periodic response of a system, the net energy transfer into or out of the system, over one period, must be zero. For the system considered here, the only nonconservative forces are the excitation and damping. As such, for the two-mode reduced order model, described by Eq.~\eqref{eq:modal_damped}, the energy relationship is
\begin{equation}
E_{D1} + E_{D2} = E_{P1} + E_{P2}\,,
\label{eq:general_energy}
\end{equation}
where $E_{Di}$ represents the energy lost by the $i^{\mbox{\scriptsize th}}$ mode due to damping, and $E_{Pi}$ represents the energy gained by the $i^{\mbox{\scriptsize th}}$ mode due to external excitation. These are computed using
\begin{equation}
E_{Di} = \int_0^T 2 \zeta_i \omega_{ni} \dot{q}_i^2 \mathrm{d}t\,,
\qquad \mbox{and} \qquad 
E_{Pi} = \int_0^T \sum_j \left[\Phi_{j,i} F_j \cos(\Omega t) \right] \dot{q}_i \mathrm{d}t\,,
\label{eq:gen_engy_components}
\end{equation}
where $T$ is the period of the response, the term in the square brackets is the $j^{\mbox{\scriptsize th}}$ external excitation force multiplied by $\Phi_{j,i}$---the mode shape of the $i^{\mbox{\scriptsize th}}$ mode at excitation point $j$. For notational convenience, $j$ will be considered equal to either 1 or 2, denoting the first and second excitation forces, for the remainder of this paper. These excitation forces are sinusoidal and at frequency $\Omega$, which is equal to the response frequency of the fundamental components of the modes (see the assumed solution in Eq.~\eqref{eq:assumed_sol}). The excitation amplitudes, $F_j$, may be positive or negative, and hence the two forces may be in- or anti-phase.

\begin{figure}[th]
\centering
\begin{tabular*}{0.95\textwidth}{@{\extracolsep{\fill}} c c}
\subfloat[]{\label{flnm1b}\includegraphics[width=0.45\textwidth]{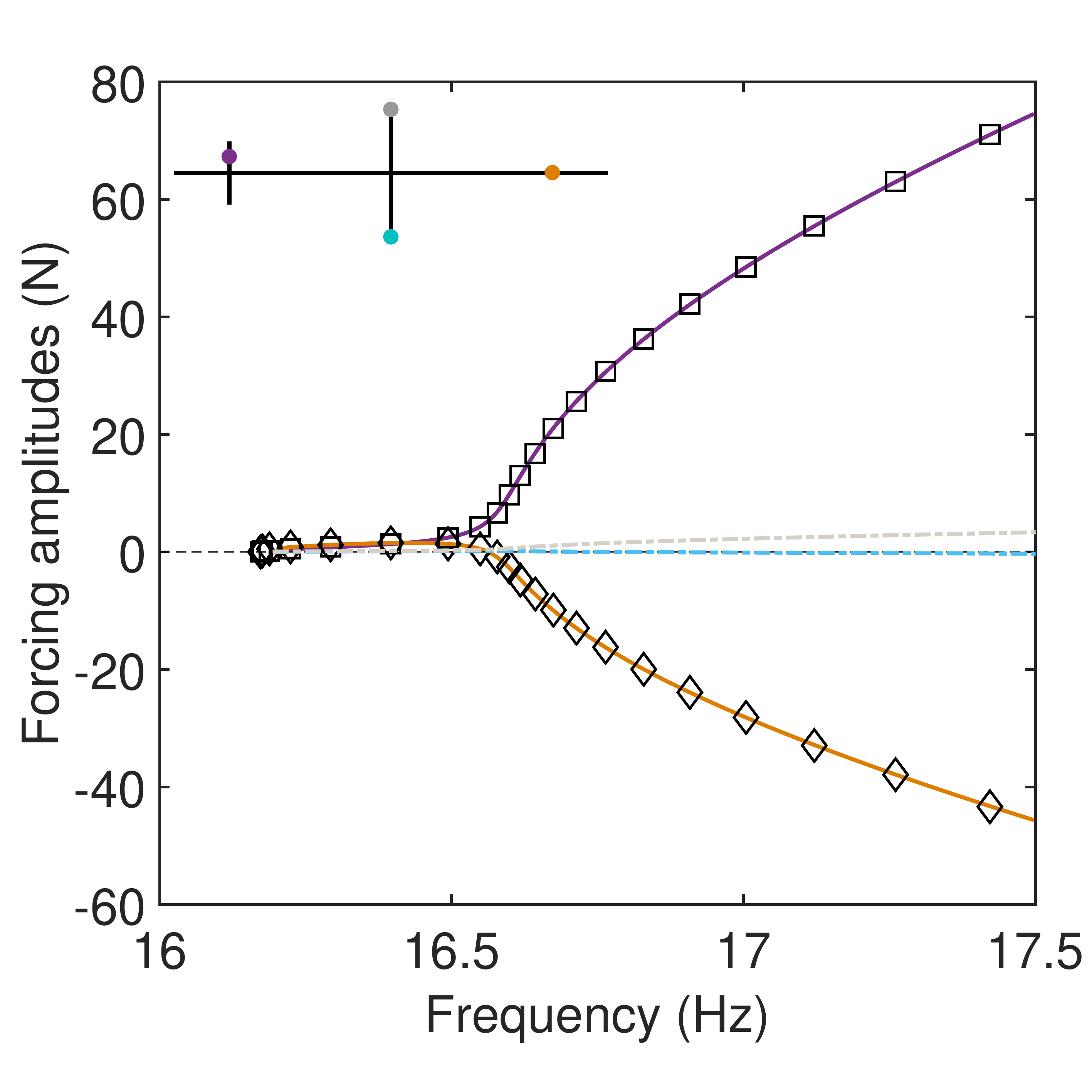}}&
\subfloat[]{\label{flnm2b}\includegraphics[width=0.45\textwidth]{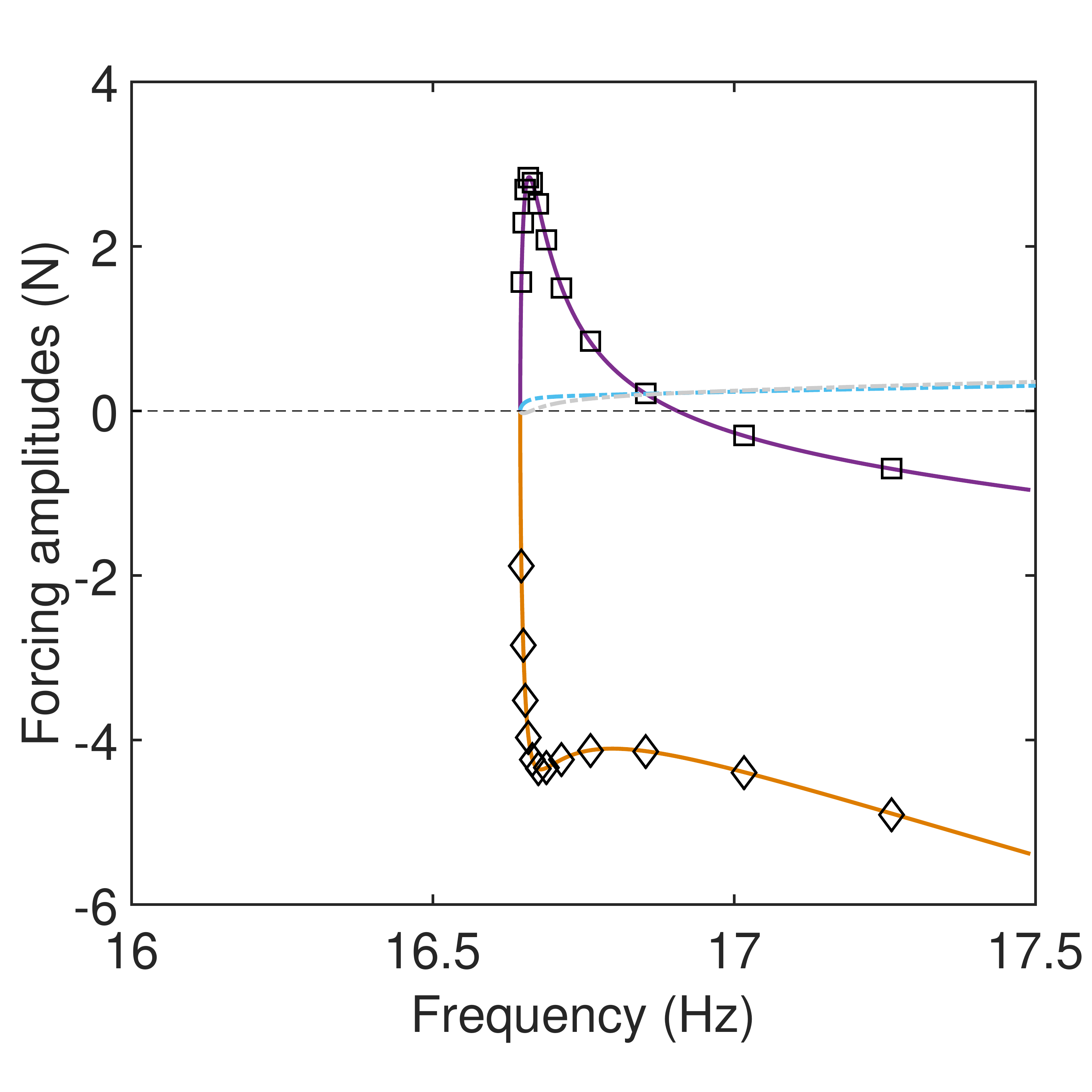}}
\end{tabular*}
\caption{Forcing amplitudes required to isolate (a) the first and (b) the second NNMs with two input forces. The first excitation considers one force on the small cross beam (\textcolor{purple}{$\boldsymbol{-}$}) and one on the main beam (\textcolor{orange}{$\boldsymbol{-}$}) as in Figure~\ref{fig:flnm12_nnms_freeratio_phys}(a). The second excitation uses forces at both tips of the cross beam (\textcolor{gray}{$\boldsymbol{-}$},\textcolor{ocean}{$\boldsymbol{-}$}). Analytical predictions for the first excitation ($\boldsymbol{\Box}$, $\boldsymbol{\Diamond}$) agree very well with numerical continuation results (solid lines).}
\label{fig:flnm12_nnms_freeratio}
\end{figure}

Whilst Eq.~\eqref{eq:general_energy} is general, and must be satisfied by any periodic response, in this paper we are specifically interested in forced responses that are equivalent to the conservative, NNM responses. As such, the additional restriction that no energy may be transferred \emph{between} the modes may be enforced. This is equivalent to stating that the energy transfer into or out of \emph{each mode} must be zero, expressed as
\begin{equation}
	E_{Di} = E_{Pi}\,,
	\qquad\mbox{for}\qquad
	i = 1,\; 2\,.
\label{eq:modal_energy_bal}
\end{equation}
Using Eqs.~\eqref{eq:gen_engy_components} and assuming two external excitation forces, Eqs.~\eqref{eq:modal_energy_bal} may be written
\begin{eqnarray}
\int_0^T 2 \zeta_i \omega_{ni} \dot{q}_i^2 \mathrm{d}t
	&\hspace{-7pt} = \hspace{-7pt}&
	\int_0^T \Phi_{1,i} F_1 \cos\left(
		\Omega t
	\right) \dot{q}_i \mathrm{d}t
	+ \int_0^T \Phi_{2,i} F_2 \cos\left(
		\Omega t
	\right) \dot{q}_i \mathrm{d}t\,,
\label{eq:energy_expr}
\end{eqnarray}
\noindent for $i=1,\;2$. The assumed solutions for the linear modes, Eq.~\eqref{eq:assumed_sol}, may now be substituted into Eqs.~\eqref{eq:energy_expr} to give
\begin{eqnarray}
\hspace{-20pt}
	2 \zeta_i \omega_{ni} \Omega^2 U_i^2 \int_0^T 
		\sin^2\left(
			\Omega t - \phi_i
		\right)
	\mathrm{d}t
	&\hspace{-7pt} = \hspace{-7pt}&
	-\left(
		\Phi_{1,i} F_1 + \Phi_{2,1} F_2 
	\right) \Omega U_i 
	\int_0^T
		\cos\left(
			\Omega t
		\right) \sin\left(
			\Omega t - \phi_i
		\right)
	\mathrm{d}t\,,
\end{eqnarray}
and, using ${ T = 2\pi\Omega^{-1} }$ these may be written
\begin{eqnarray}
\hspace{-20pt}
	2 \zeta_i  \omega_{ni} \Omega U_i^2
	&\hspace{-7pt} = \hspace{-7pt}&
	\left(
		\Phi_{1,i} F_1 + \Phi_{2,i} F_2 
	\right) U_i \sin\left( \phi_i \right).
\label{eq:energy_rel_sin}
\end{eqnarray}
It is assumed that the forced response is in quadrature with the responses, i.e.~${ \phi_i = \pm\pi/2 }$. As such, ${ \sin\left(\phi_i\right) = p_i = \pm 1 }$. Furthermore, it has been shown (Section~\ref{sec:11beam}) that, on the NNMs, the linear modes are in- or anti-phase, denoted ${ p = +1 }$ and ${ p = -1 }$ respectively. As such, we may define ${ p_2 = p p_1 }$, and hence Eqs.~\eqref{eq:energy_rel_sin} may be written
\begin{subequations}
\begin{eqnarray}
\hspace{-20pt}
	2 \zeta_1  \omega_{n1} \Omega U_1
	&\hspace{-7pt} = \hspace{-7pt}&
	\left(
		\Phi_{1,1} F_1 + \Phi_{2,1} F_2 
	\right) p_1\,,
\label{eq:energy_rel_sin_1_p1}
\\
\hspace{-20pt}
	2 \zeta_2  \omega_{n2} \Omega U_2
	&\hspace{-7pt} = \hspace{-7pt}&
	\left(
		\Phi_{1,2} F_1 + \Phi_{2,2} F_2
	\right) p p_1\,.
\label{eq:energy_rel_sin_2_p1}
\end{eqnarray}%
\label{eq:energy_rel_sin_p1}%
\end{subequations}%
As the left-hand sides of Eqs.~\eqref{eq:energy_rel_sin_p1} are always positive, it follows that the right-hand sides must also be positive. As $F_1$ and $F_2$ may be positive or negative, this highlights that a choice may be made regarding the sign of $p_1$ --- i.e.~a positive $F_1$ and $F_2$ with ${ p_1 = +1 }$ is equivalent to a negative $F_1$ and $F_2$ with ${ p_1 = -1 }$. Here, we choose to define ${ p_1 = +1 }$ such that Eqs.~\eqref{eq:energy_rel_sin_p1} are written
\begin{equation}
	\left[\begin{array}{cc}
		\Phi_{1,1} 	& \Phi_{1,2}
	\\	\Phi_{2,1}	& \Phi_{2,2}
	\end{array}\right]^{*}
	\left(\begin{array}{c}
		F_1
	\\	F_2
	\end{array}\right)
	= \left(\begin{array}{c}
		P_1
	\\	P_2
	\end{array}\right) = \left(\begin{array}{c}
		2 \zeta_1  \omega_{n1} \Omega U_1
	\\	2 p \zeta_2  \omega_{n2} \Omega U_2
	\end{array}\right),
\label{eq:energy_rel}%
\end{equation}%
where $(.)^*$ denotes the transpose operator and $P_1$ and $P_2$ are the modal forcing amplitudes. Eq.~\eqref{eq:energy_rel} is similar to the one obtained for linear systems~\cite{GeradinBook}, except that the amplitudes $U_i$ are now a solution of the nonlinear Eqs.~\eqref{eq:resonant_sols}. Substituting such a solution into Eq.~\eqref{eq:energy_rel}, the forcing amplitudes that are required to obtain that response may be computed. To validate analytical predictions, Figure~\ref{fig:flnm12_nnms_freeratio} compares the force amplitudes predicted by Eqs.~\eqref{eq:energy_rel} ($\boldsymbol{\Box}$, $\boldsymbol{\Diamond}$) and the results obtained using numerical continuation for the case of the excitation considered in Figure~\ref{fig:flnm12_nnms_freeratio_phys}, revealing a very good agreement. 

Note that not all pairs of physical coordinates represent suitable excitation locations. In particular, when both excitation locations lead to near-parallel modal force vectors, the modal forcing required to isolate some NNM motions may be reached only with unrealistically large excitation amplitudes. As already observed, one of the input forces may also change sign (and hence phase) to compensate for the other input forces and maintain the appropriate excitation. As shown by Eqs.~\eqref{eq:energy_rel}, this depends on the choice of excitation locations (and hence the coefficients $\Phi_{j,i}$) which can be selected to minimise this effect. This is, for instance, the case of the blue and grey locations at the tips of the main cross beam (Figure~\ref{fig:flnm12_nnms_freeratio}(a)). However, the placement of an exciter at these locations would be considered impractical given the large displacements and rotations that the structure exhibits at these points. 

In Section~\ref{sec:pred_energy}, the energy balance approach introduced in this section will be revisited for the case of a single input force and used to analyse the accuracy of quadrature curves.

\section{Phase quadrature---restricted excitation conditions}\label{sec:nnm_1fphys}
The previous section has shown that the NNMs of a system can be captured accurately using a sufficient number of suitably-positioned input forces. Furthermore, the number of excitation forces does not have to be equal to the number of DOFs in the system, only to the number of strongly-interacting modes. However, in many applications, structures may have a high modal density leading to a large number of interacting modes, whilst the number of possible excitation locations is normally limited (often to one) for experimental convenience and practicality. As such, this section investigates the case where fewer excitation forces than interacting modes are considered.

\begin{figure}[th]
\begin{center}
\includegraphics[width=1.0\textwidth]{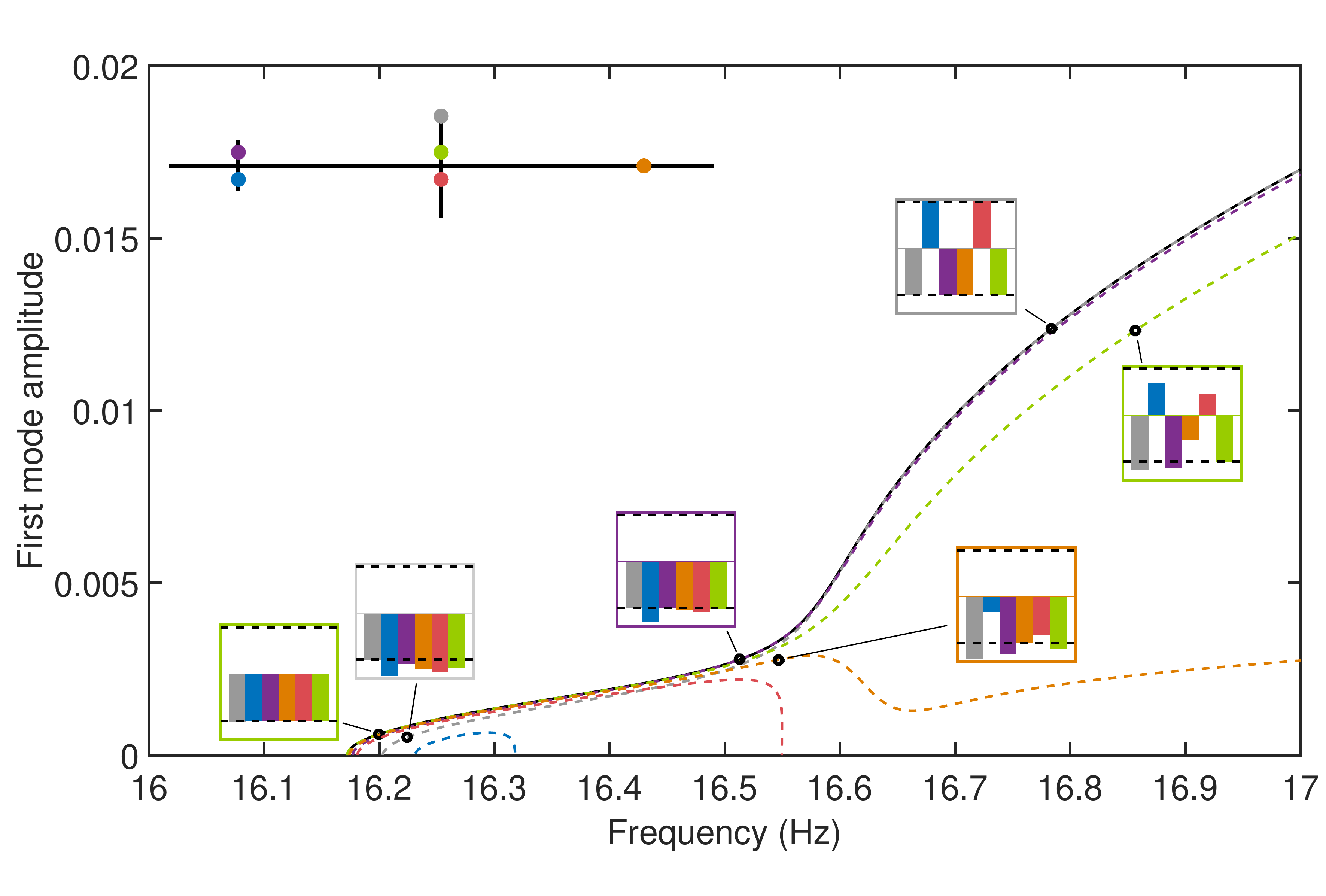}
\caption{Comparison between the first NNM (\textcolor{black}{$\boldsymbol{-}$}) and the quadrature curves obtained with a single-harmonic, single-point excitation applied at the locations reported on the structure schematic. Insets show the phase of the response at all the excitation locations for the periodic solutions highlighted by $\boldsymbol{\circ}$. Quadrature curves, structure locations and insets are coloured accordingly.}
\label{fig:11beam_nnm1_1fphys}
\end{center}
\end{figure}

A single-harmonic, single-point excitation force is now applied to the beam structure. The quadrature condition imposes a phase difference of $\pm 90$ deg. between the phase of the excitation and the phase of the response of the structure at the excited DOF. The quadrature curves obtained when considering six different excitation locations are represented by dashed lines in Figure~\ref{fig:11beam_nnm1_1fphys}. Each quadrature curve is coloured according to the location of the excitation, as shown in the structure schematic inset. The first NNM of the underlying conservative system is represented in solid black (\textcolor{black}{$\boldsymbol{-}$}).

The six excitation locations lead to six different quadrature curves --- all approximating the same NNM but with various degrees of accuracy, depending on the response frequency and amplitude. For instance, the red and orange quadrature curves in Figure~\ref{fig:11beam_nnm1_1fphys} show a completely different trend from the NNM as frequency increases, despite a good agreement at first. Similar observations can be made of the green curve, which captures the NNM qualitatively but fails to capture it quantitatively, despite the location of the excitation providing an almost perfect isolation of the first linear mode ($P_1 = 0.9999$ and $P_2 = 6\times 10^{-3}$). Conversely, the grey curve fails to capture the low-amplitude region (including natural frequency) of the NNM but becomes accurate at higher amplitudes. This variation in the high-amplitude results regardless of the low-amplitude accuracy shows that linear modal analysis considerations alone cannot be considered to appropriately select the position of the input force, at least for the coupled system considered here.

Inspection of the quadrature curves alone does not allow us to identify which curve resembles the NNM most closely. As with linear systems~\cite{Wright99}, an inaccurate isolation of the NNM will result in phase differences in the response across the structure and hence the quadrature condition will not be met for all coordinates. This observation is illustrated in Figure~\ref{fig:11beam_nnm1_1fphys} with bar charts showing the phases of the response of the structure at the different excitation locations. The colour of the outer box surrounding the bar charts matches the quadrature curve under consideration. Figure~\ref{fig:11beam_nnm1_1fphys} shows that the phase differences observed across the structure cannot generally be overlooked as they are symptomatic of potentially large errors between the quadrature and the NNM. In an experiment, they should be carefully monitored. Section~\ref{sec:pred_energy} will combine this notion of phase difference with energy balance arguments to further analyse the source of (in)accurate quadrature and discuss the selection of suitable excitation locations. Note that the response co-located with the excitation will always appear at $\pm90$ deg., as their phase difference was used to define the quadrature condition that is tracked using numerical continuation.

\begin{figure}[th]
\centering
\includegraphics[width=0.95\textwidth]{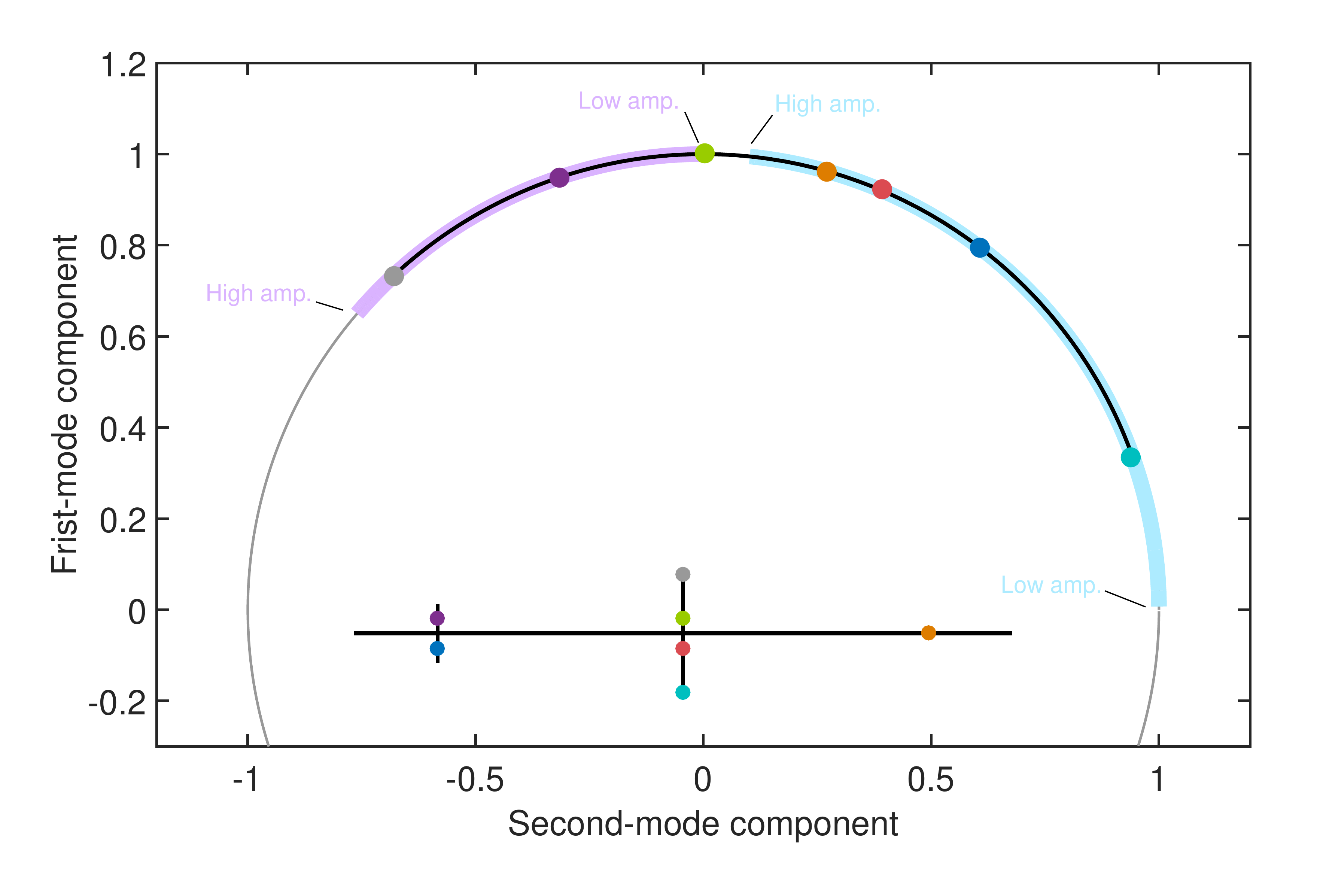}
\caption{Modal input force ratios required to obtain the first (\textcolor{pink}{$\hrectangleblack$}) and second (\textcolor{lblue}{$\hrectangleblack$}) NNMs up to a frequency of 20 Hz. Coloured bullet points show the force obtained for the excitation locations reported on the structure schematic. (\textcolor{black}{$\boldsymbol{-}$}) All possible modal force ratios achievable with an excitation applied to only one vertical DOF of the finite element model.}
\label{fig:11beam_fphys_circle}
\end{figure}

From Figure~\ref{fig:11beam_nnm1_1fphys}, it is clear that the location of the excitation plays a key role in the resulting NNM prediction. In fact, it dictates the ratio of input forces applied to modes 1 and 2. Figure~\ref{fig:11beam_fphys_circle} shows the ratio of modal input forces required, theoretically, to obtain the first (\textcolor{pink}{$\hrectangleblack$}) and second (\textcolor{lblue}{$\hrectangleblack$}) NNMs up to a frequency of 20 Hz, as calculated using Eqs.~\eqref{eq:energy_rel}. Modal force vectors were normalised for presentation on the unit circle. The force vectors imposed by the single excitation locations considered in Figure~\ref{fig:11beam_nnm1_1fphys} are also represented as bullet points coloured according to their location.

Each location provides an excitation that is best suited to only one NNM motion. Similar observations were previously made in Figure~\ref{fig:flnm12_nnms_freeratio_phys}(b) where the excitation applied to the small cross was the only one required to isolate the NNM motion at 16.55 Hz (see $\boldsymbol{+}$). At other frequencies, forcing at both excitation points were needed. More generally, excitation locations also lead to either positive (\textcolor{green}{$\boldsymbol{\bullet}$},\textcolor{purple}{$\boldsymbol{\bullet}$},\textcolor{gray}{$\boldsymbol{\bullet}$}) or negative (\textcolor{ocean}{$\boldsymbol{\bullet}$},\textcolor{blue}{$\boldsymbol{\bullet}$},\textcolor{myred}{$\boldsymbol{\bullet}$},\textcolor{orange}{$\boldsymbol{\bullet}$}) ratios between modal forces, and hence are more suited to either the first or second NNM. While Figure~\ref{fig:11beam_fphys_circle} shows the exact ratio of first and second modes needed to capture particular NNM motions, it does not mean that a given force ratio cannot be used to approximate the NNM elsewhere.

The evolution of the exact modal forcing follows the evolution of the right-hand side of Eqs.~\eqref{eq:energy_rel} and can therefore be understood a priori through the evolution of the mode shapes (Figure~\ref{fig:nnms}). For instance, as the response of the first NNM evolves from bending to torsion, the location of the ideal excitation evolves towards locations providing larger torsional components. Note, perhaps surprisingly, that opposite locations such as blue (\textcolor{blue}{$\boldsymbol{\bullet}$}) and purple (\textcolor{purple}{$\boldsymbol{\bullet}$}), or red  (\textcolor{myred}{$\boldsymbol{\bullet}$}) and green (\textcolor{green}{$\boldsymbol{\bullet}$}), lead to significantly different modal forces and hence quadrature curves. This comes from the geometric asymmetry of the linear modes that is induced by the asymmetric position of the concentrated masses. The geometric asymmetry in the system is slight and would therefore be difficult to detect in practice. This suggest that the appropriate excitation of more-complicated structure could be very challenging. 

Figure~\ref{fig:11beam_fphys_circle} also includes a solid black line which represents all modal force vectors realisable by applying the excitation force on one of the vertical DOFs of the finite element model. Interestingly, some modal excitations cannot be achieved with a single input force and would require the superposition of two forces. For instance, there is no location across the structure where a vertical excitation would provide sufficient torsion and small-enough bending components to isolate the second NNM at low amplitude. It also appears that positioning actuators on the best locations could be challenging in practice; in particular when the considered locations involve large displacements/rotations as, for instance, on the cross beam (\textcolor{gray}{$\boldsymbol{\bullet}$},\textcolor{green}{$\boldsymbol{\bullet}$},\textcolor{myred}{$\boldsymbol{\bullet}$},\textcolor{ocean}{$\boldsymbol{\bullet}$}). Such practical constraints further reduce the range of forcing vectors available for NNM appropriation.

\begin{figure}[th]
\centering
\includegraphics[width=1.0\textwidth]{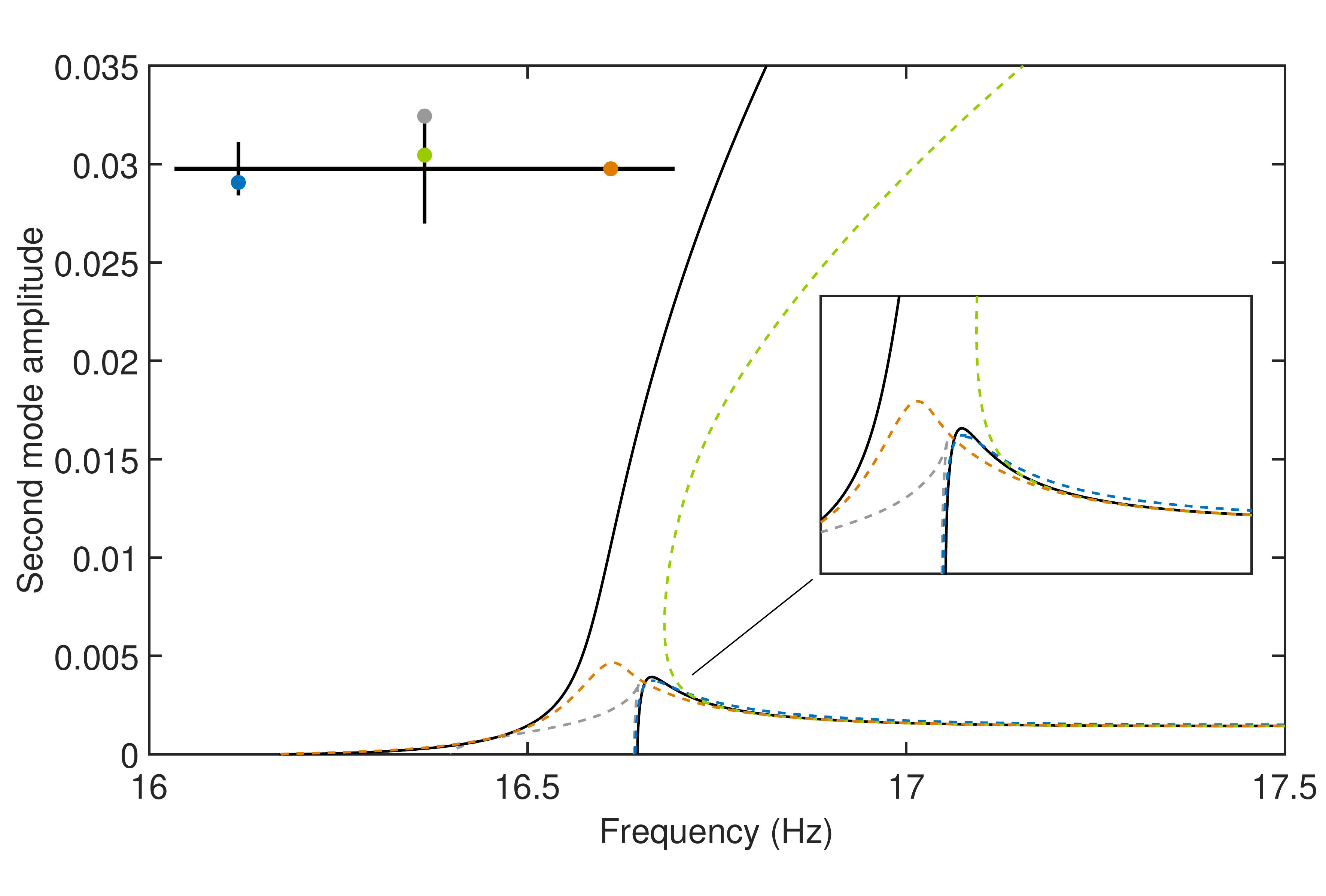}
\caption{Comparison between the second NNM (\textcolor{black}{$\boldsymbol{-}$}) and quadrature curves obtained with a single-harmonic, single-point excitation applied at the locations reported on the structure schematic. Note that the first NNM is also shown for comparison.}
\label{fig:11beam_nnm2_1fphys}
\end{figure}

Figure~\ref{fig:11beam_nnm2_1fphys} shows four quadrature curves obtained in the vicinity of the second NNM. The orange quadrature curve, which was previously presented in Figure~\ref{fig:11beam_nnm1_1fphys}, is shown to join the second NNM when its mode shape includes sufficient bending component. The trend of the grey quadrature curve is qualitatively different from the NNM. The blue location provides an excellent estimate of the entire NNM, although small inaccuracies at low amplitude  are noticeable. Similar inaccuracies are observable for all excitation locations as the theoretically required excitation cannot be reached with a single input force (Figure~\ref{fig:11beam_fphys_circle}). 

In addition to the quadrature curves found at low amplitude, other curves that exist only for sufficiently large response amplitudes were found. This is, for instance, the case of the green quadrature curve shown in Figure~\ref{fig:11beam_nnm2_1fphys}. The computation of that curve was initiated at high amplitude on the second NNM, i.e. where the modal force vector is close to the modal force vector given by an excitation at the green location. Half of this curve corresponds to the second NNM whereas the other half does not resemble either of the NNMs. The presence of such isolated quadrature curves --- disconnected from any linear quadrature condition --- can be attributed to the significant change in the NNM mode shapes with amplitude (Figure~\ref{fig:nnms}).

\section{Influence of energy transfer on quadrature accuracy}\label{sec:pred_energy}
Section~\ref{sec:nnm_1fphys} has shown that the quadrature curves obtained with a single-force excitation can be significantly different from the NNMs of the underlying conservative system. At best, quadrature curves were shown to capture the NNM in a limited range of amplitudes. Furthermore, it was shown that linear considerations may not allow for the identification of suitable excitation locations, in particular for high-amplitude responses where the deformation of the structure has evolved from the underlying linear mode shapes.

This section further develops the energy arguments presented in Section~\ref{sec:theory}\ref{sec:balance} to demonstrate that energy transfer between modes can be used to understand the accuracy of the quadrature condition. First, the amplitude of the single-force excitation required to obtain a specific NNM motion is calculated (Section~\ref{sec:pred_energy}\ref{sec:balance2}), then the phase difference between modal responses introduced by energy transfer into and out of a modal coordinate is derived analytically (Section~\ref{sec:pred_energy}\ref{sec:phase_error}). Analytical results are compared to the numerical simulation results of Section~\ref{sec:nnm_1fphys}.

\subsection{Analytical prediction of the single-force excitation amplitude}\label{sec:balance2}
The energy balancing approach, previously considered in Section \ref{sec:theory}\ref{sec:balance} for the case where two excitation forces are provided, is now applied for the case of a single excitation point. The relationship between the forcing and damping energy, defined in Eqs.~\eqref{eq:general_energy} and \eqref{eq:gen_engy_components} enforce that there is no net energy gain, or loss, from the system during one period of motion. For the case of a single force, the energy input in each mode cannot be modified separately and Eq.~\eqref{eq:general_energy} cannot be separated into Eqs.~\eqref{eq:modal_energy_bal}. As such, a single input force cannot control the transfer of energy from one mode to the other, which may therefore be non-zero. For an NNM response, it is required that the net energy transfer between the modes is zero; hence a single-point forcing cannot always reach an NNM precisely, as observed in previous sections.

Combining Eqs.~\eqref{eq:general_energy} and \eqref{eq:gen_engy_components} for the case of a single forcing gives
\begin{equation}
	  \int_0^T 2 \zeta_1 \omega_{n1} \dot{q}_1^2 \mathrm{d}t
	+ \int_0^T 2 \zeta_2 \omega_{n2} \dot{q}_2^2 \mathrm{d}t
	= \int_0^T  \Phi_{1,1} F_1 \cos(\Omega t) \dot{q}_1 \mathrm{d}t
	+ \int_0^T  \Phi_{1,2} F_1 \cos(\Omega t) \dot{q}_2 \mathrm{d}t\,.
\label{eq:energy_expr_sum}
\end{equation}
The assumed solutions for the linear modes, Eq.~\eqref{eq:assumed_sol}, are substituted into Eq.~\eqref{eq:energy_expr_sum} to give
\begin{eqnarray}
&&\hspace{-40pt}
	2 \zeta_1 \omega_{n1} \Omega^2 U_1^2 \int_0^T 
		\sin^2\left(
			\Omega t - \phi_1
		\right)
	\mathrm{d}t
	+ 2 \zeta_2 \omega_{n2} \Omega^2 U_2^2 \int_0^T
		\sin^2\left(
			\Omega t - \phi_2
		\right)
	\mathrm{d}t
	= 
\notag\\
&&\hspace{-20pt}
	-\Phi_{1,1} F_1 \Omega U_1 \int_0^T
		\cos\left(
			\Omega t
		\right) \sin\left(
			\Omega t - \phi_1
		\right)
	\mathrm{d}t
	- \Phi_{1,2} F_1 \Omega U_2 \int_0^T
		\cos\left(
			\Omega t
		\right) \sin\left(
			\Omega t - \phi_2
		\right)
	\mathrm{d}t,
\end{eqnarray}
and, using 
${ T = 2\pi\Omega^{-1} }$ this is written
\begin{equation}
	 2 \zeta_1 \omega_{n1} \Omega U_1^2
	+ 2 \zeta_2 \omega_{n2} \Omega U_2^2
	= \Phi_{1,1} F_1 U_1\sin\left(\phi_1\right)
	+ \Phi_{1,2} F_1 U_2\sin\left(\phi_2\right).
\label{eq:final_energy_w/phase}
\end{equation}

\begin{figure}[th]
\centering
\begin{tabular*}{0.95\textwidth}{@{\extracolsep{\fill}} c c}
\subfloat[]{\label{f1a}\includegraphics[width=0.45\textwidth]{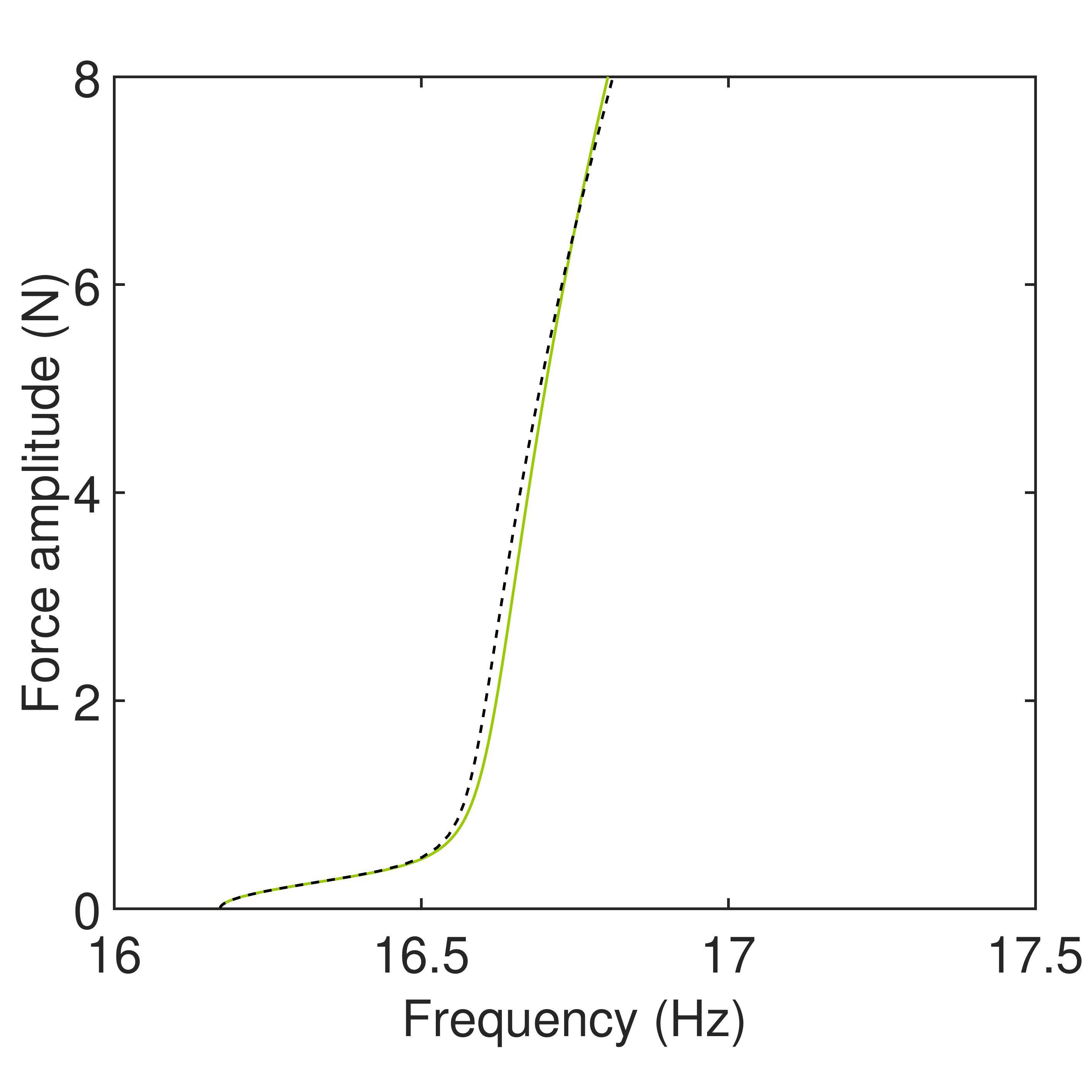}}&
\subfloat[]{\label{f1b}\includegraphics[width=0.45\textwidth]{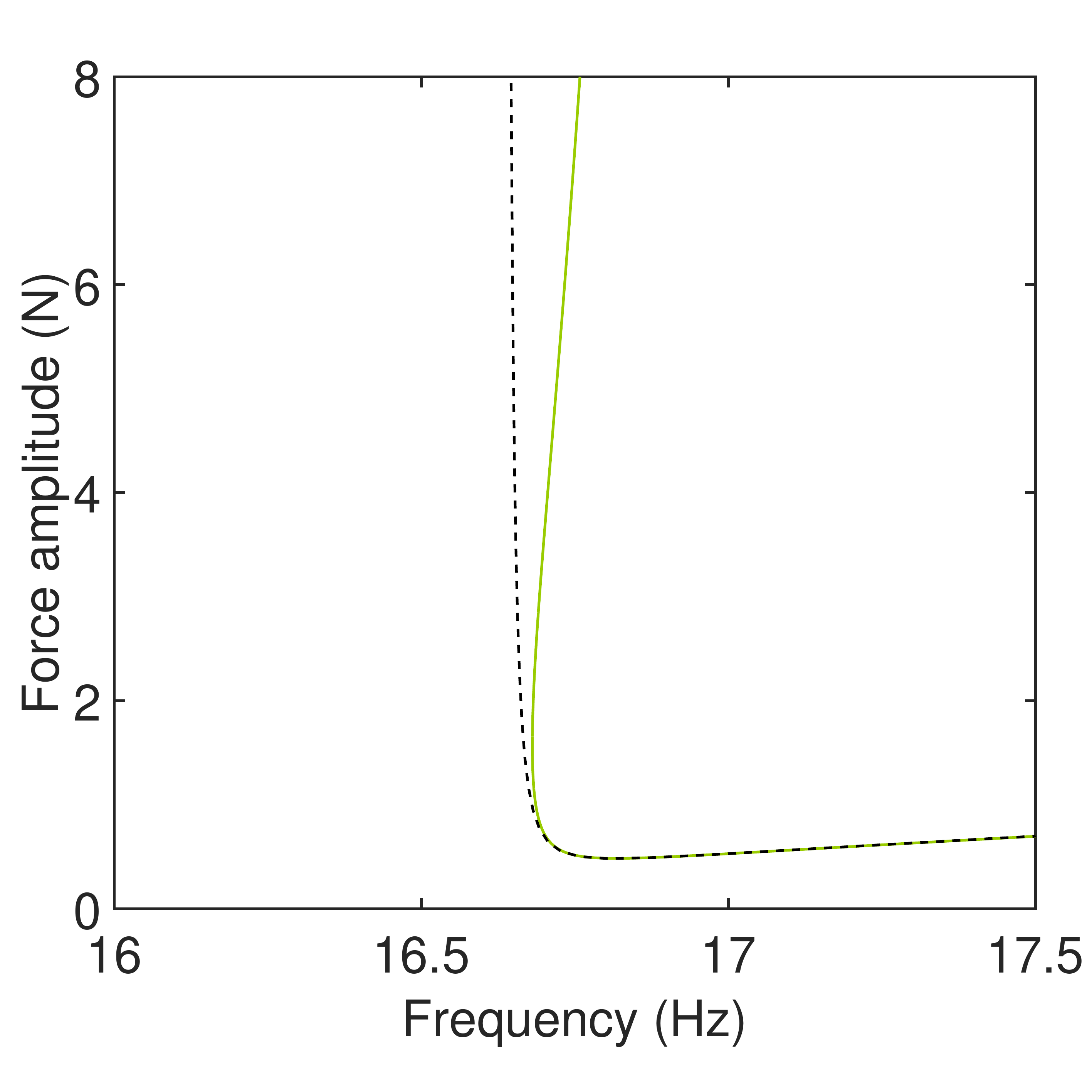}}
\end{tabular*}
\caption{Force amplitude required to isolate the first (a) and second (b) NNMs of the clamped-clamped cross beam structure in the case of a single-point, single-harmonic excitation applied to the green location (\textcolor{green}{$\boldsymbol{\bullet}$}). Analytical predictions given by Eq.~\eqref{eq:F1} (\textcolor{black}{$\boldsymbol{--}$}) and results obtained using numerical continuation (\textcolor{green}{$\boldsymbol{-}$}).} 
\label{fig:F1}
\end{figure}

As in Section~\ref{sec:theory}\ref{sec:balance}, it is assumed that the forced response is in quadrature with the applied excitation, i.e.~${ \phi_i = \pm\pi/2 }$, and that ${ \sin\left(\phi_1\right) = p_1 = +1 }$ and ${ \sin\left(\phi_2\right) = p p_1 }$. Hence Eq.~\eqref{eq:final_energy_w/phase} can be written as
\begin{equation}
	 2 \zeta_1 \omega_{n1} \Omega U_1^2
	+ 2 \zeta_2 \omega_{n2} \Omega U_2^2
	= p_1 F_1\left(
		  \Phi_{1,1} U_1 
		+ p \Phi_{1,2} U_2
	\right),
\label{eq:final_energy_w/p}
\end{equation}
which may also be rearranged to give
\begin{equation}
	F_1 = 
	\dfrac{
		  2 \Omega\left(
		  	  \zeta_1 \omega_{n1} U_1^2
			+ \zeta_2 \omega_{n2} U_2^2
		\right)
	}{
		  \Phi_{1,1} U_1 
		+ p \Phi_{1,2} U_2
	}.
\label{eq:F1}
\end{equation}
As in the case of two excitation forces, the solutions of Eqs.~\eqref{eq:resonant_sols} may be substituted into Eq.~\eqref{eq:F1} to find the required forcing amplitude. Figure~\ref{fig:F1} compares the force amplitude calculated using Eq~\eqref{eq:F1} (\textcolor{black}{$\boldsymbol{--}$}) with the force amplitude obtained using numerical continuation (\textcolor{green}{$\boldsymbol{-}$}) for the case of the green excitation location. A very good agreement between analytical predictions and numerical results is found in regions where the quadrature curve matches the NNM as shown in Figures~\ref{fig:F1}(a) and~\ref{fig:F1}(b) below 16.5 Hz and above 16.7 Hz, respectively. Errors appear in regions where the quadrature curve deviates from the NNM. This is because both modal coordinates are no longer in quadrature with the excitation, which violates one of the assumptions leading to Eq.~\eqref{eq:F1}.

\subsection{Phase error due to energy transfer}\label{sec:phase_error}
As discussed in~\cite{Hill15,Ehrhardt17sub}, energy transfer between modes results in phase differences between modes that are different to that on the NNMs, i.e.~${ \phi_1 - \phi_2 \neq \pm \pi }$. Such phase differences were already noticed in Figure~\ref{fig:11beam_nnm1_1fphys} where all the considered coordinates were not oscillating synchronously in quadrature with the excitation. 

To analytically predict the phase difference between modes, we consider the different possible mechanisms for transferring energy into and out of a mode, i.e. the modal forcing, the modal damping and the nonlinear couplings with the other modes. For any periodic response, regardless of its proximity to an NNM, the net energy transfer for any individual mode must be zero over one period. Considering the first mode of the cross-beam, whose coupling terms are defined in Eq.~\eqref{eq:app_ROM_BB}, the condition for zero net energy transfer may be written
\begin{equation}
	\int_0^T
	  	\left[
	  		2 \zeta_1 \omega_{n1} \dot{q}_1
	  	\right] \dot{q}_1
	\mathrm{d}t
	+ \int_0^T
		\left[
		  	3\gamma_2 q_1^2 q_2
			+ \gamma_3 q_1 q_2^2
			+ \gamma_4 q_2^3
		\right] \dot{q}_1
	\mathrm{d}t
	= \int_0^T
		\left[
	  		\Phi_{1,1} F_1 \cos(\Omega t)
	  	\right] \dot{q}_1
	\mathrm{d}t.
\label{eq:q1_energy_expr}
\end{equation}
Substituting the assumed solutions for the linear modes, Eq.~\eqref{eq:assumed_sol}, into Eq.~\eqref{eq:q1_energy_expr} and evaluating the integrals (for
${ T = 2\pi\Omega^{-1} }$) gives
\begin{equation}
	  2 \zeta_1 \omega_1 U_1 \Omega
	+ \dfrac{3}{4} \gamma_2 U_1^2 U_2 \sin\left(\hat{\phi}_d\right)
	+ \dfrac{1}{4} \gamma_3 U_1 U_2^2 \sin\left(2\hat{\phi}_d\right)
	+ \dfrac{3}{4} \gamma_4 U_2^3 \sin\left(\hat{\phi}_d\right)
	= \Phi_{1,1} F_1 \sin\left(\phi_1\right),
\label{eq:ET_expr}
\end{equation}
where
${ \phi_1 - \phi_2 = \hat{\phi}_d }$ has been used. Here $\hat{\phi}_d$ denotes the \emph{phase-error}---i.e.~the phase difference that results from the energy transfer between the modes. Now, using the previous convention that ${ \sin(\phi_1) = +1 }$, Eq.~\eqref{eq:ET_expr} may be written
\begin{equation}
	\sin(\hat{\phi}_d)
	= \dfrac{
		4\left[
			\Phi_{1,1} F_1 - 2 \zeta_1 \omega_1 U_1 \Omega
		\right]
	}{
		\left[
			  3 \gamma_2 U_1^2 
			+ 2 \gamma_3 U_1 U_2 \cos(\hat{\phi}_d)
			+ 3 \gamma_4 U_2^2
		\right] U_2
	}.
    \label{eq:phase_diff}
\end{equation}
Note that using the second equation of motion would lead to an expression similar to~\eqref{eq:phase_diff}, providing an identical estimate of the phase error, $\hat{\phi}_d$.

Eq.~\eqref{eq:phase_diff} provides an estimation of the phase difference between modes given a particular NNM motion, the location, and the amplitude of the excitation and the damping properties of the system. Using the analytic approximation of the NNM (Eqs.~\eqref{eq:resonant_sols}) and the estimation of the input force required to isolate a NNM motion (Eq.~\eqref{eq:F1}), the nonlinear Eq.~\eqref{eq:phase_diff} may be solved for $\hat{\phi}_d$. This phase-error represents the accuracy of the quadrature in representing a NNM --- large values of $\hat{\phi}_d$ correspond to less accurate quadrature conditions. Note that, in some cases, no real solutions for $\hat{\phi}_d$ exist, representing a phase-error with a magnitude exceeding $\pi/2$---i.e.~where the quadrature solution is highly inaccurate and far from the NNM it represents.
\begin{figure}
\centering
\subfloat[]{\label{nfbb1}\includegraphics[width=0.9\textwidth]{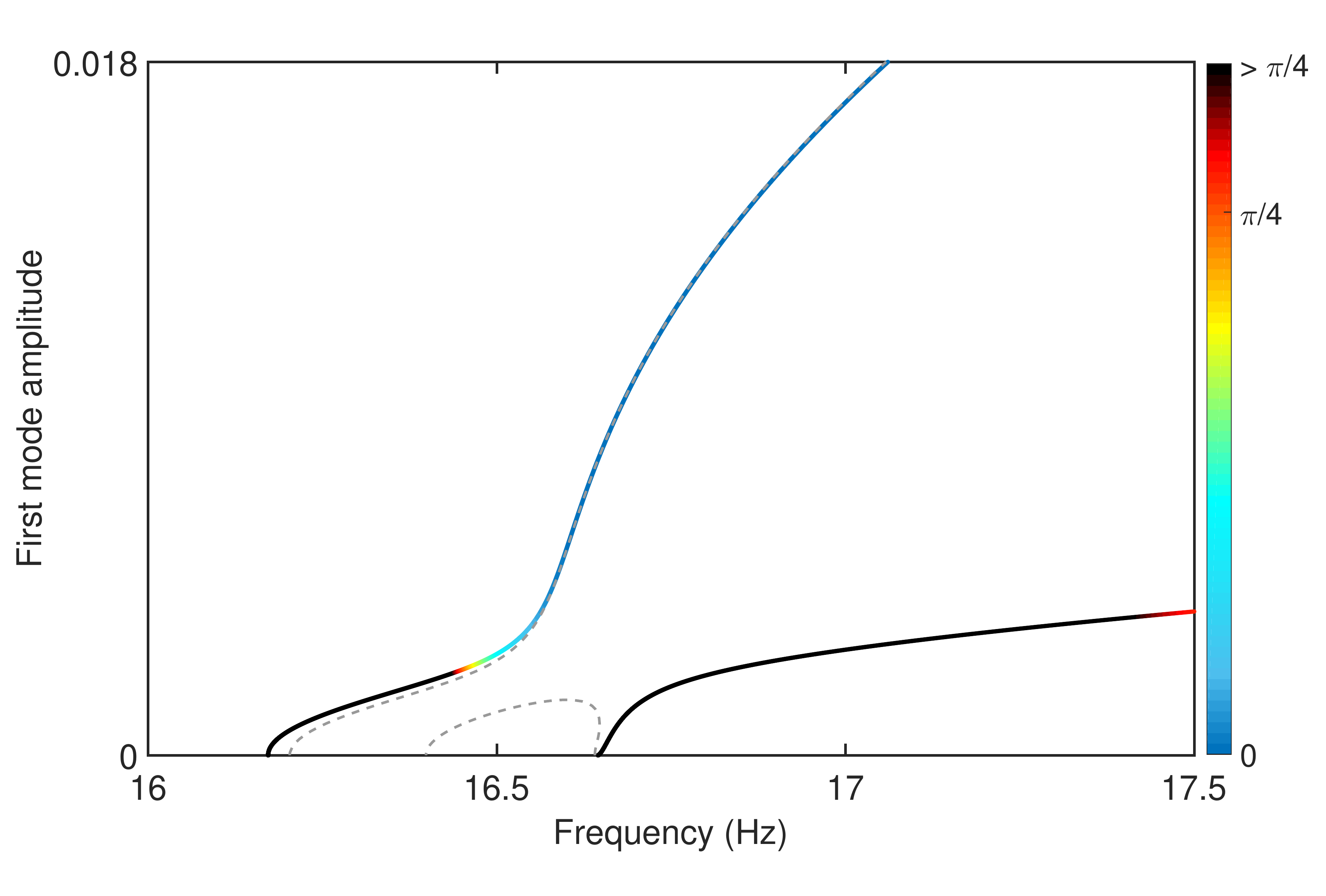}}\\
\subfloat[]{\label{nfbb4}\includegraphics[width=0.9\textwidth]{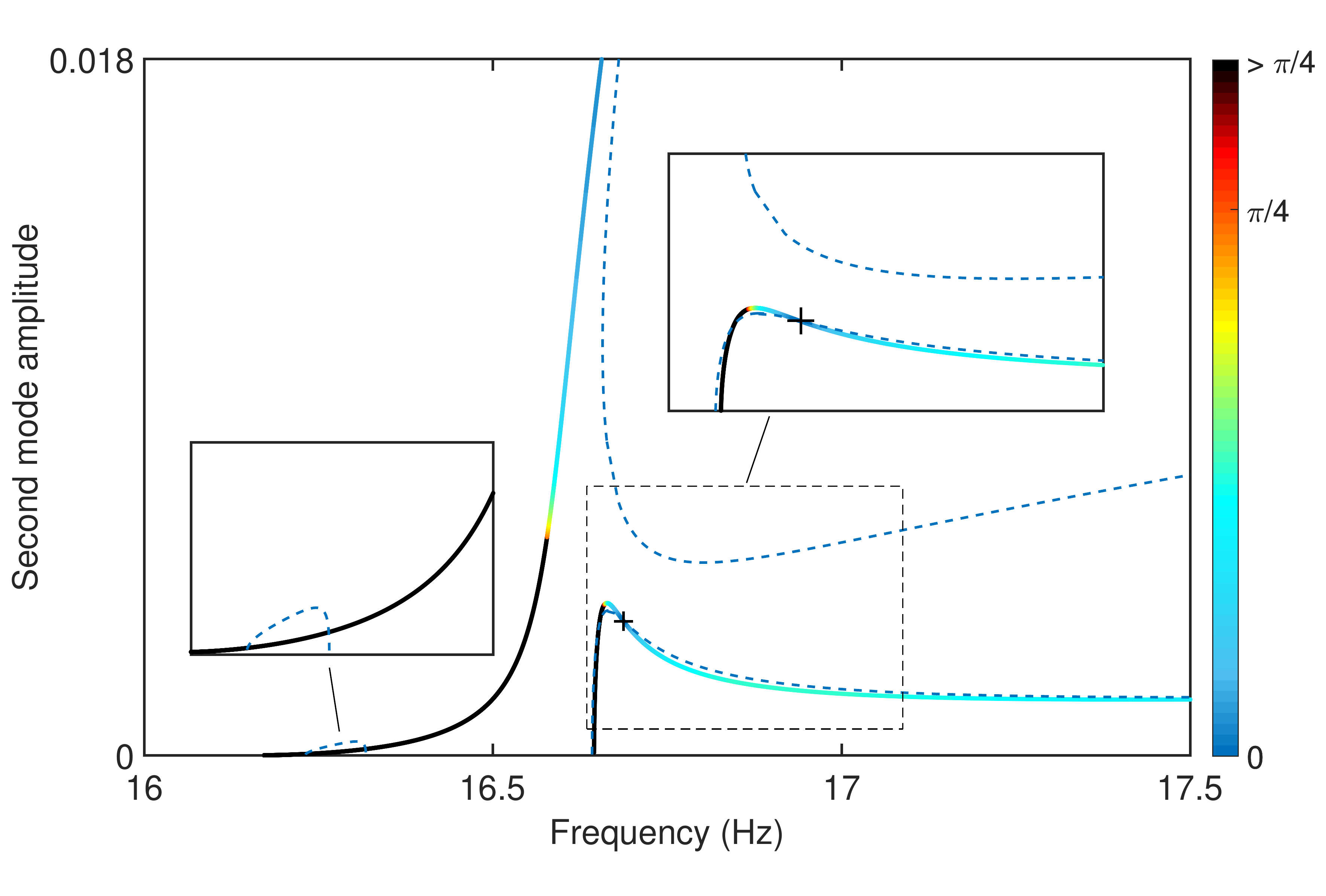}}
\caption{Regions with large phase errors predict where quadrature curves will be inaccurate. NNMs (solid lines) are coloured according to the magnitude of the phase error. All phase errors larger than $\pi/4$ are coloured in black. Quadrature curves (dashed lines) are coloured according to the excitation location in Figure~\ref{fig:11beam_fphys_circle}: (a) main cross beam (\textcolor{gray}{$\boldsymbol{\bullet}$}), (b) small cross beam (\textcolor{blue}{$\boldsymbol{\bullet}$}). The intersection between the quadrature curve and the second NNM is at 16.69 Hz ($\boldsymbol{+}$).}
\label{fig:NFBBs}
\end{figure}

The NNMs of the beam structure are represented in Figure~\ref{fig:NFBBs} with a colour scale showing the magnitude of the phase error $\hat{\phi}_d$. Overall, the analytical predictions correctly identify the regions where the quadrature curves obtained for a single-force excitation are (in)accurate. For instance, Figure~\ref{fig:NFBBs}(a) shows that an excitation applied to the grey location (\textcolor{gray}{$\boldsymbol{\bullet}$}) is inadequate for the beginning of the first NNM and for the whole of the second NNM. Conversely, the ability to recover a good estimation of the first NNM at high amplitudes is also well captured. 

The blue location is also shown to be inadequate for low-amplitude responses (Figure~\ref{fig:NFBBs}(b)). For the second NNM, a small phase error is predicted where the quadrature curve intersects the actual NNM ($\boldsymbol{+}$). This intersection also occurs for the projection in the first mode. Beyond that region, moderately small errors are predicted. Small phase errors are also predicted at high-amplitude on the first NNM. This observation was exploited to find the isolated quadrature condition and track its evolution with forcing amplitude using numerical continuation (the upper dashed curve \textcolor{blue}{$\boldsymbol{--}$}).
\begin{figure}
\centering
\includegraphics[width=0.95\textwidth]{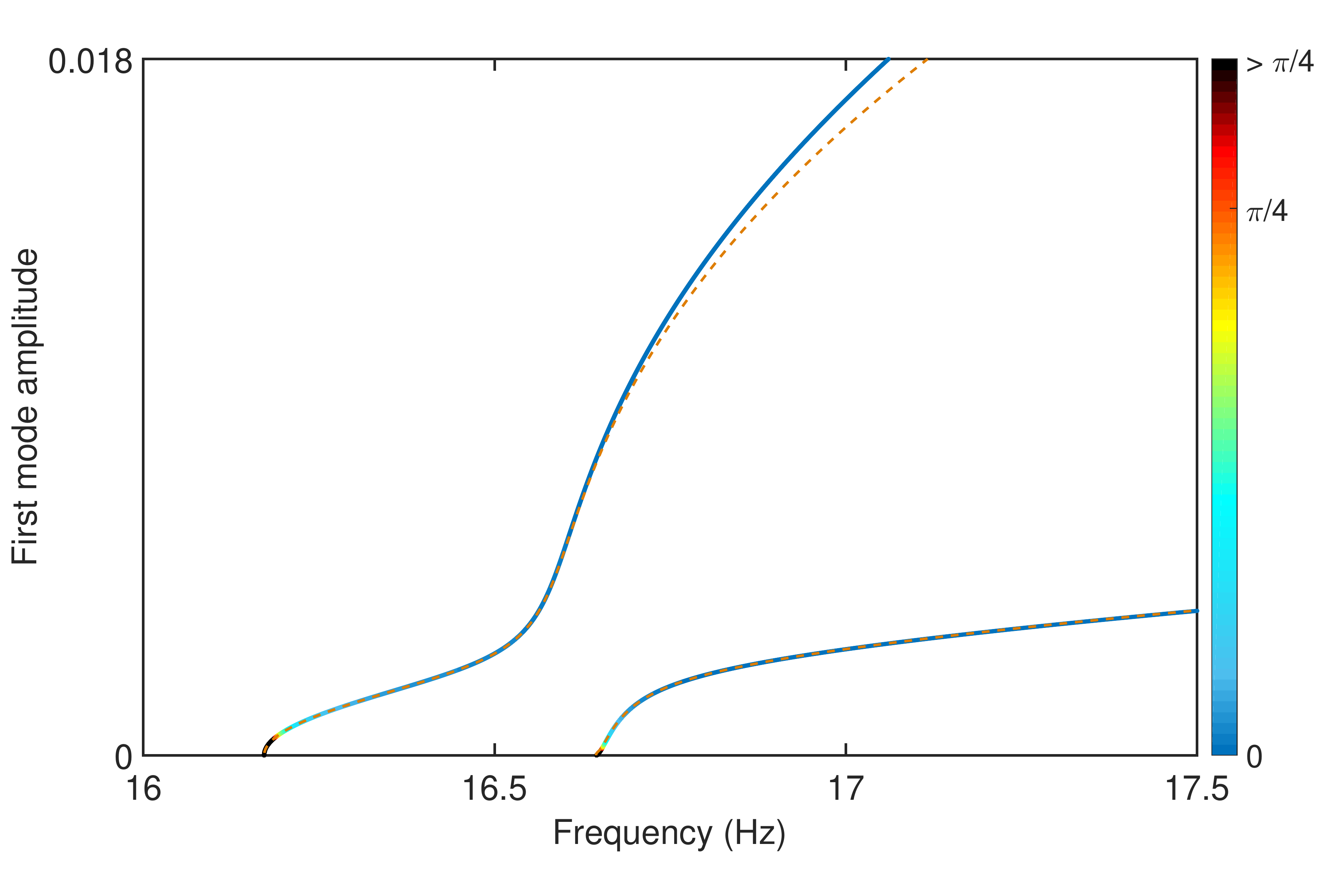}
\caption{Comparison between the NNMs of the structure (coloured according to the magnitude of the phase error) and the quadrature curves numerically obtained when damping ratios are ten times smaller. Single-point harmonic excitation applied vertically to the orange location (\textcolor{orange}{$\boldsymbol{\bullet}$}). All phase errors larger than $\pi/4$ are coloured in black.}
\label{fig:NFBBs_C10}
\end{figure}

Eq.~\eqref{eq:phase_diff} shows that the phase error depends on the damping properties of the system. As such, excitation locations leading to inaccurate quadrature curves for some damping values may lead to more accurate curves for lower levels of damping. This is illustrated in Figure~\ref{fig:NFBBs_C10} where the accuracy of the curve obtained with a single-point excitation at the orange location is significantly improved when the damping in the system is ten times smaller than that considered in Table~\ref{tab:11beam_all_prop}. More specifically, compared to Figures~\ref{fig:11beam_nnm1_1fphys} and~\ref{fig:11beam_nnm2_1fphys}, the transition from the first to the second NNM no longer exists and one of the quadrature curves matches a larger potion of the first NNM. The second NNM is also very well captured by a second quadrature curve that starts at low response amplitude in the vicinity of the second natural frequency. Quadrature conditions in that region were not found in Figures~\ref{fig:11beam_nnm1_1fphys} and~\ref{fig:11beam_nnm2_1fphys}. Note that the damping characteristics reported in Table~\ref{tab:11beam_all_prop} are reasonable for engineering systems as obtained from experimental tests on a similar structure~\cite{Ehrhardt17sub,RensonISMA2016}.

\subsection{Comparison between different excitation locations}
\begin{figure}
\centering
\includegraphics[width=0.95\textwidth]{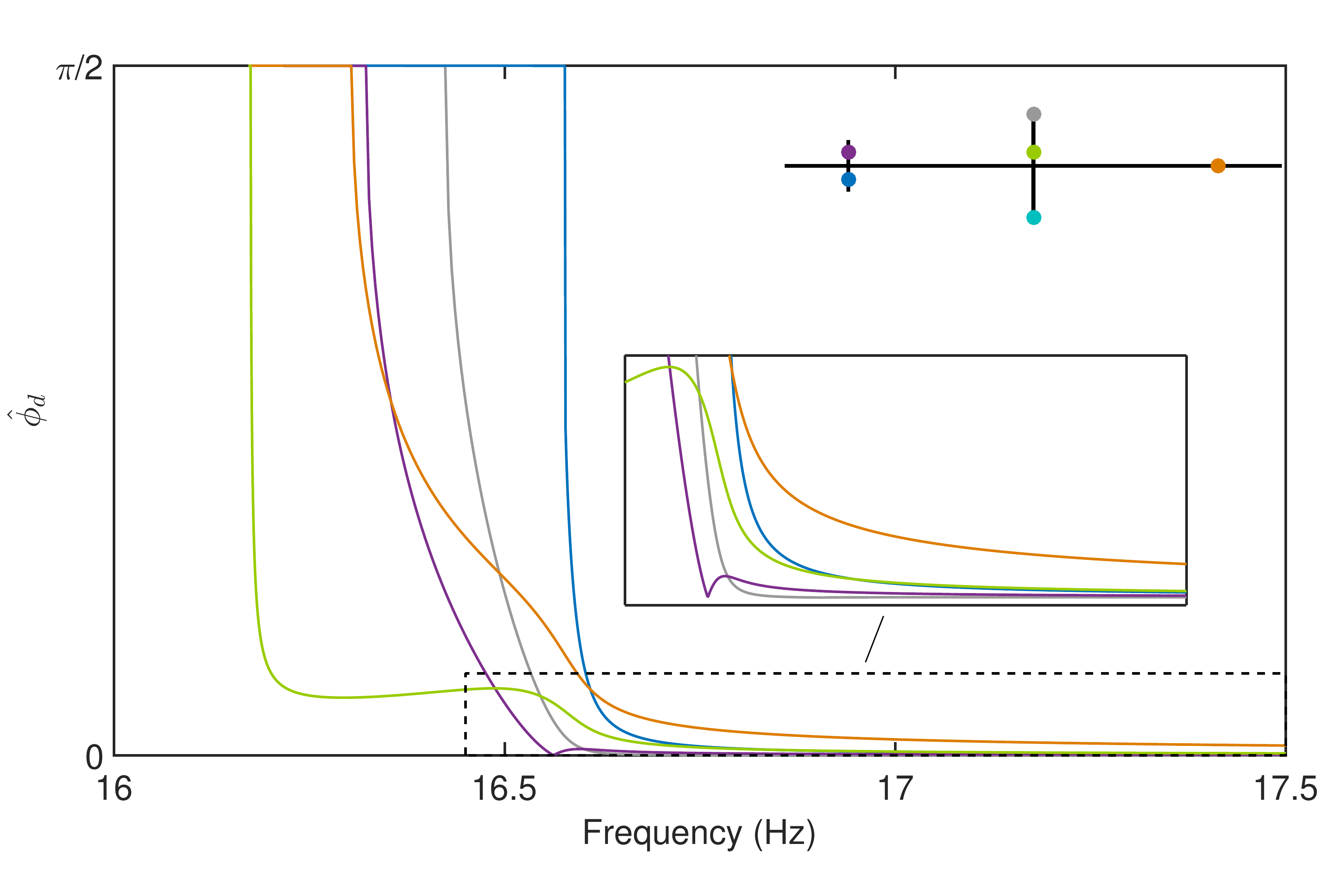}
\caption{Evolution of the phase error $\hat{\phi}_d$, Eq.~\eqref{eq:phase_diff}, for the first NNM and the different force locations reported on the structure schematic.}
\label{fig:comp_phiD}
\end{figure}

Eq.~\eqref{eq:phase_diff} allows us to qualitatively estimate the accuracy of the quadrature curve without directly computing it. Compared to numerical continuation which requires the calculation of the entire quadrature curves, the analytic approach enables a rapid comparison of different excitation locations across the structure. This is shown in Figure~\ref{fig:comp_phiD} for the excitation locations considered in the structure schematic inset in the figure. The results are consistent with the previous analysis of these quadrature curves. Close to the first natural frequency, the green location provides the best quadrature condition, followed by the purple and then the grey locations as frequency increases. It can also be observed that the blue location can provide a better quadrature curve than the orange and green locations for frequencies larger than about 17 Hz. This observation is also consistent with the presence of the isolated quadrature curve shown in Figure~\ref{fig:NFBBs}(b). 

\begin{figure}[t]
\setlength{\belowcaptionskip}{-10pt}
\centering
\includegraphics[width=1.0\textwidth]{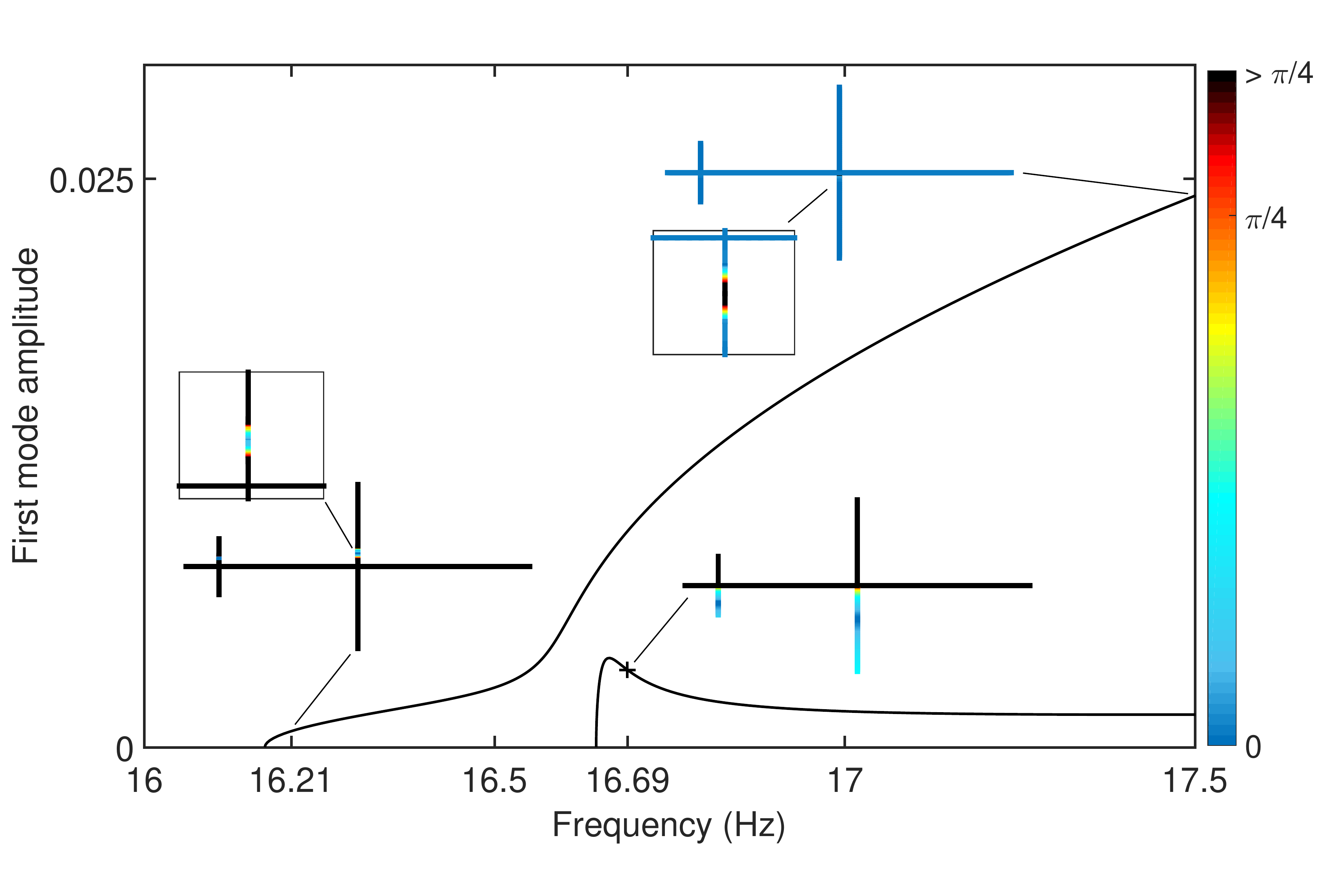}
\caption{Phase error $\hat{\phi}_d$ obtained when considering each node of the finite-element model as the excitation location. Three particular NNM motions are considered. Each node in the associated structure schematic is coloured according to the magnitude of the phase error.}
\label{fig:comp_phiD_stru}
\end{figure}

In Figure~\ref{fig:comp_phiD_stru}, each node of the structure is coloured according to the phase error obtained when applying the input force on the vertical DOF associated with that node. Three different NNM motions are considered. For the first NNM, at 16.61 Hz, only the locations in the vicinity of the green excitation point (see, for instance, Figure~\ref{fig:11beam_nnm1_1fphys}) provide relatively-low phase errors, which is consistent with earlier observations. At 17.5 Hz, the phase error is small almost everywhere on the structure, which is expected given the behaviour of the phase error as amplitude increases (Figure~\ref{fig:comp_phiD}). A small region on the cross beam shows, however, high phase error values. This is because the denominator in Eq.~\eqref{eq:F1} is very small. So the force $F_1$ required to isolate the NNM motion considered is large and hence the numerator of Eq.~\eqref{eq:phase_diff} is also large. The NNM motion taken at 16.69 Hz on the second NNM is the same as the one highlighted in Figure~\ref{fig:NFBBs}(b)($\boldsymbol{+}$). One of the locations leading to a small phase error corresponds to the blue location, for which the quadrature curve was shown to coincide with the second NNM (Figure~\ref{fig:NFBBs}(b)).

The successful analytical predictions provided by Eq.~\eqref{eq:phase_diff} shows that energy transfer between modes is a key factor in the accuracy of quadrature curves. However, Eq.~\eqref{eq:phase_diff} may not be directly considered as a quantitative estimation of the quadrature accuracy. Take, for example, the predictions made for low amplitudes in Figure~\ref{fig:NFBBs_C10}. The quadrature curve is regarded as inaccurate whereas the frequency-amplitude evolution of the NNM is well reproduced. At high amplitude, the phase error remains very small whereas the quadrature curve clearly starts to deviate from the NNM.

These issues arise from the nature of the coupling between the modes. In fact, all the curves in Figure~\ref{fig:comp_phiD} exhibit high (low) phase errors for low (high) response amplitudes. For near-linear responses, energy transfer due to modal coupling are weaker and so a comparatively small energy transfer may only be achieved by a large change in the phase between the modes. Conversely, a high-amplitude response is likely to exhibit a strong modal coupling, and thus may transfer a very large amount of energy, leading to smaller phase errors.

The above issues could also be partly attributed to the projection used in Figure~\ref{fig:NFBBs_C10}. At low response amplitude, the relative error appears more clearly in the second modal coordinate (not represented for conciseness). However, the contribution of the second mode to the response of the system is at least two orders of magnitude smaller than the first mode, such that the frequency-amplitude dependence of the NNM remains well captured and the phase-error criterion can be considered as overly pessimistic. At high response amplitude, modal coordinates contribute equally to the response such that the phase-error criterion is overly optimistic.

\section{Conclusion}\label{sec:conclu}
Quadrature curves extracted from the forced and damped response of a structure are often considered to represent the nonlinear normal modes of that structure. This paper has shown that significant errors between quadrature curves and nonlinear normal modes can exist when a limited number of input forces is used and, in particular, when the excitation is restricted to a single input force as it is frequently the case in practice. The structure analysed --- a doubly-clamped cross beam structure --- was particularly challenging due to the presence of strong modal couplings (1:1 interaction), the significant evolution of the nonlinear mode shapes with amplitude and the spatial decoupling between bending and torsion introduced by the cross shape of the physical structure. As such, quadrature curves obtained with a single force were found to be very sensitive to the location of the excitation.

By considering the energy transferred between the degrees of freedom of the system, we were able to construct an analytical expression to estimate the phase error introduced by particular forcing conditions. Although this phase error may not directly map to the accuracy of the quadrature curves, it proved to be a good accuracy indicator. The phase error was exploited to identify regions of the nonlinear normal modes that can be well captured with given excitation conditions. It was also used to highlight the forcing locations that could be used to capture a particular nonlinear normal mode motion.

\section*{Data Statement}
The data presented in this work are openly available from the University of Bristol repository at https://doi.org/10.5523/bris.2ryns8g7ircfv2dsoguq44f1pd.

\section*{Acknowledgement}
L.R. has received funding from the Royal Academy of Engineering (research fellowship RF1516/15/11), D.A.E. is funded by the Engineering Nonlinearity EPSRC Program Grant EP/K003836/1, D.A.W.B. is funded by the EPSRC grant EP/K032738/1, S.A.N. by the EPSRC fellowship EP/K005375/1. We gratefully acknowledge the financial support of the Royal Academy of Engineering and EPSRC.

\appendix
\section{Analytical approximation of NNMs}\label{app:NNMs}
The equations of motion describing the backbone curves of the reduced order model are found by combining Eqs.~\eqref{eq:modal_damped} and~\eqref{eq:nq_rom} and removing the forcing and damping terms to give
\begin{subequations}
\begin{eqnarray}
	\ddot{q}_1
	+ \omega_{n1}^2 q_1
	+ \gamma_1 q_1^3
	+ 3\gamma_2 q_1^2 q_2
	+ \gamma_3 q_1 q_2^2
	+ \gamma_4 q_2^3
	&\hspace{-6pt}=\hspace{-6pt}& 0,
\\
	\ddot{q}_2
	+ \omega_{n2}^2 q_2
	+ \gamma_2 q_1^3 
	+ \gamma_3 q_1^2 q_2 
	+ 3 \gamma_4 q_1 q_2^2 
	+ \gamma_5 q_2^3
	&\hspace{-6pt}=\hspace{-6pt}& 0.
\end{eqnarray}%
\label{eq:app_ROM_BB}%
\end{subequations}%
It is then assumed that the responses of the modes may be written
\begin{equation}
	q_i
	\approx u_i
	= U_i \cos\left(
		\Omega t - \phi_i
	\right),
\label{eq:app_assumed_sol}
\end{equation}
where $U_i$ and $\phi_i$ denote the fundamental response amplitude and phase of the $i^{\mbox{\scriptsize th}}$ mode respectively, and $\Omega$ represents the response frequency. These assumed responses are now substituted into Eqs.~\eqref{eq:app_ROM_BB}. Then, using a harmonic balance approximation, the non-resonant terms (i.e.~those terms that do not resonate at frequency $\Omega$) are removed, giving
\begin{subequations}
\begin{eqnarray}
	  \ddot{u}_1
	+ \omega_{n1}^2 u_1
	+ \dfrac{U_1}{4}\left[
		3 \gamma_1 U_1^2 + 2 \gamma_3 U_2^2
	\right] \cos\left(
		\Omega t - \phi_1
	\right)
	+ \dfrac{3 U_2}{4} \left[
		2 \gamma_2 U_1^2 + \gamma_4 U_2^2
	\right] \cos\left(
		\Omega t - \phi_2
	\right)
	+ \qquad &&
\\
	\dfrac{3 \gamma_2}{4} U_1^2 U_2 \cos\left(
		\Omega t - 2\phi_1 + \phi_2
	\right)
	+ \dfrac{\gamma_3}{4} U_1 U_2^2 \cos\left(
		\Omega t + \phi_1 - 2\phi_2
	\right)
	= 0, &&
\notag
\\
	\ddot{u}_2
	+ \omega_{n2}^2 u_2
	+ \dfrac{3 U_1}{4} \left[
		\gamma_2 U_1^2 + 2 \gamma_4 U_2^2
	\right] \cos\left(
		\Omega t - \phi_1
	\right)
	+ \dfrac{U_2}{4} \left[
		2 \gamma_3 U_1^2 + 3 \gamma_5 U_2^2
	\right] \cos\left(
		\Omega t - \phi_2
	\right)
	+ \qquad &&
\\
	\dfrac{\gamma_3}{4} U_1^2 U_2 \cos\left(
		\Omega t - 2\phi_1 + \phi_2
	\right)
	+ \dfrac{3 \gamma_4}{4} U_1 U_2^2 \cos\left(
		\Omega t + \phi_1 - 2\phi_2
	\right)
	= 0. &&
\notag
\end{eqnarray}%
\label{eq:app_resonant_EOM}%
\end{subequations}%
Note that the quadratic nonlinear terms originally present in Eqs.~\eqref{eq:modal_damped} have been ignored in Eqs.~\eqref{eq:app_ROM_BB} because they would not appear in the resonant equations of motion~\cite{TouzeCISM}. Quadratic terms are also ignored in the energy analysis as they have a negligible influence on the results and removing them considerably simplifies the analytical expressions.

Using the assumed solutions, Eqs.~\eqref{eq:app_resonant_EOM} may then be written
\begin{subequations}
\begin{eqnarray}
&&\hspace{-20pt}
	\left\{
		\left(
			\omega_{n1}^2 - \Omega^2
		\right) U_1
		+ \dfrac{3 \gamma_1}{4} U_1^3
		+ \dfrac{3 U_2}{4} \left[
			3 \gamma_2 U_1^2
			+ \gamma_4 U_2^2
		\right] \cos\left( \phi_d \right)	
		+ \dfrac{\gamma_3}{4} U_1 U_2^2 \left[
			2 + \cos\left(2\phi_d\right)
		\right]			
	\right\} \cos\left(
		\Omega t - \phi_1
	\right)
\notag\\
&&\hspace{50pt}
	- \dfrac{U_2}{4} \left\{
		3 \left[
		  	\gamma_2 U_1^2
			+ \gamma_4 U_2^2
		\right] \sin\left(\phi_d\right)		
		+ \gamma_3 U_1 U_2 \sin\left( 2\phi_d\right)		
	\right\} \sin\left(
		\Omega t - \phi_1
	\right)	
	= 0,
\\
&&\hspace{-20pt}
	\left\{
		\left(
			\omega_{n2}^2 - \Omega^2
		\right) U_2
		+ \dfrac{3 U_1}{4} \left[
			\gamma_2 U_1^2
			+ 3 \gamma_4 U_2^2
		\right] \cos\left(\phi_d\right)
		+ \dfrac{\gamma_3}{4} U_1^2 U_2 \left[
			2 + \cos\left( 2\phi_d \right)
		\right]
		+ \dfrac{3 \gamma_5}{4} U_2^3		
	\right\} \cos\left(
		\Omega t - \phi_2
	\right)
\notag\\
&&\hspace{50pt}
	+ \dfrac{U_1}{4} \left\{
		  3 \left[
		  	\gamma_2 U_1^2
			+ \gamma_4 U_2^2
		\right] \sin\left(\phi_d\right)
		+ \gamma_3 U_1 U_2 \sin\left( 2\phi_d \right)		
	\right\} \sin\left(
		\Omega t - \phi_2
	\right)
	= 0.
\end{eqnarray}%
\label{eq:app_resonant_EOM_2}%
\end{subequations}%
where
${ \phi_d = \phi_1 - \phi_2 }$.
Noting that the terms in the curly braces, i.e.~$\{\bullet\}$, corresponding to 
$\sin\left(\Omega t - \phi_1\right)$ and
$\sin\left(\Omega t - \phi_2\right)$ are equal, Eqs.~\eqref{eq:app_resonant_EOM_2} are satisfied when
\begin{subequations}
\begin{eqnarray}
	\left(
		\omega_{n1}^2 - \Omega^2
	\right) U_1
	+ \dfrac{3 \gamma_1}{4} U_1^3
	+ \dfrac{3 U_2}{4} \left[
		3 \gamma_2 U_1^2
		+ \gamma_4 U_2^2
	\right] \cos\left( \phi_d \right)	
	+ \dfrac{\gamma_3}{4} U_1 U_2^2 \left[
		2 + \cos\left(2\phi_d\right)
	\right]
	&\hspace{-6pt}=\hspace{-6pt}& 0,\qquad
\label{eq:app_resonant_comp_c1}%
\\
	\left(
		\omega_{n2}^2 - \Omega^2
	\right) U_2
	+ \dfrac{3 U_1}{4} \left[
		\gamma_2 U_1^2
		+ 3 \gamma_4 U_2^2
	\right] \cos\left(\phi_d\right)
	+ \dfrac{\gamma_3}{4} U_1^2 U_2 \left[
		2 + \cos\left( 2\phi_d \right)
	\right]
	+ \dfrac{3 \gamma_5}{4} U_2^3
	&\hspace{-6pt}=\hspace{-6pt}& 0,\qquad
\label{eq:app_resonant_comp_c2}%
\\
	3 \left[
	  	\gamma_2 U_1^2
		+ \gamma_4 U_2^2
	\right] \sin\left(\phi_d\right)		
	+ \gamma_3 U_1 U_2 \sin\left( 2\phi_d\right)		
	&\hspace{-6pt}=\hspace{-6pt}& 0.\qquad
\label{eq:app_resonant_comp_s}%
\end{eqnarray}%
\label{eq:app_resonant_comp}%
\end{subequations}%
Equation.~\eqref{eq:app_resonant_comp_s} is satisfied when
${ \phi_d = 0 }$, i.e.~the case where the two modes are in-phase, or when ${ \phi_d = \pi }$, where the two modes are in anti-phase. Substituting these solutions into Eqs.~\eqref{eq:app_resonant_comp_c1} and~\eqref{eq:app_resonant_comp_c2} gives
\begin{subequations}
\begin{eqnarray}
	\left(
		\omega_{n1}^2 - \Omega^2
	\right) U_1
	+ \dfrac{3 \gamma_1}{4} U_1^3
	+ p \dfrac{3 U_2}{4} \left[
		3 \gamma_2 U_1^2
		+ \gamma_4 U_2^2
	\right]
	+ \dfrac{3 \gamma_3}{4} U_1 U_2^2
	&\hspace{-6pt}=\hspace{-6pt}& 0,
\label{eq:app_resonant_sol1}%
\\
	\left(
		\omega_{n2}^2 - \Omega^2
	\right) U_2
	+ p \dfrac{3 U_1}{4} \left[
		\gamma_2 U_1^2
		+ 3 \gamma_4 U_2^2
	\right] 
	+ \dfrac{3 \gamma_3}{4} U_1^2 U_2
	+ \dfrac{3 \gamma_5}{4} U_2^3
	&\hspace{-6pt}=\hspace{-6pt}& 0,
\label{eq:app_resonant_sol2}%
\end{eqnarray}%
\label{eq:app_resonant_sols}%
\end{subequations}%
where ${ p = +1 }$ when the modes are in-phase (second NNM) and ${ p = -1 }$ when the modes are in anti-phase (first NNM). Equations~\eqref{eq:app_resonant_sols} may now be solved to find the fundamental amplitudes of the modes, $U_1$ and $U_2$, in terms of the response frequency, $\Omega$.

\bibliographystyle{unsrt}
\bibliography{mybib}

\begin{thebibliography}{10}

\bibitem{Noel14}
J.~P. No\"el, L.~Renson, and G.~Kerschen.
\newblock Complex dynamics of a nonlinear aerospace structure: Experimental
  identification and modal interactions.
\newblock {\em Journal of Sound and Vibration}, 333(12):2588--2607, 2014.

\bibitem{VakakisTET}
A.~F. Vakakis, O.~Gendelman, L.~A. Bergman, D.~M. McFarland, G.~Kerschen, and
  Y.~S. Lee.
\newblock {\em Nonlinear Targeted Energy Transfer in Mechanical and Structural
  Systems}, volume 156 of {\em Solid Mechanics and Its Applications}.
\newblock Springer Netherlands, 2009.

\bibitem{Qalandar14}
K.~R. Qalandar, B.~S. Strachan, B.~Gibson, M.~Sharma, A.~Ma, S.~W. Shaw, and
  K.~L. Turner.
\newblock Frequency division using a micromechanical resonance cascade.
\newblock {\em Applied Physics Letters}, 105(24):244103, 2014.

\bibitem{Rhoads08}
J.F. Rhoads, S.W. Shaw, and K.L. Turner.
\newblock Nonlinear dynamics and its applications in micro- and nanoresonators.
\newblock In {\em Proceedings of the 2008 ASME Dynamic Systems and Control
  Conference}, Ann Arbor, Michigan, USA, 2008.

\bibitem{Renson15}
L.~Renson, J.~P. No\"el, and G.~Kerschen.
\newblock Complex dynamics of a nonlinear aerospace structure: numerical
  continuation and normal modes.
\newblock {\em Nonlinear Dynamics}, 79(2):1293--1309, 2015.

\bibitem{Ducceschi13}
M.~Ducceschi, C.~Touz\'e, S.~Bilbao, and C.~J. Webb.
\newblock Nonlinear dynamics of rectangular plates: investigation of modal
  interaction in free and forced vibrations.
\newblock {\em Acta Mechanica}, pages 1--20, 2013.

\bibitem{Liu16}
X.~Liu, A.~Cammarano, D.~J. Wagg, S.~A. Neild, and R.~J. Barthorpe.
\newblock {$N-1$} modal interactions of a three-degree-of-freedom system with
  cubic elastic nonlinearities.
\newblock {\em Nonlinear Dynamics}, 83(1):497--511, 2016.

\bibitem{Kerschen09}
G.~Kerschen, M.~Peeters, J.~C. Golinval, and A.~F. Vakakis.
\newblock Nonlinear normal modes, part {I}: A useful framework for the
  structural dynamicist.
\newblock {\em Mechanical Systems and Signal Processing}, 23(1):170--194, 2009.

\bibitem{GeradinBook}
M.~G\'eradin and D.~Rixen.
\newblock {\em Mechanical Vibrations: Theory and application to structural
  dynamics}.
\newblock John Wiley, 1997.

\bibitem{Vakakis97}
A.~F. Vakakis.
\newblock Non-linear normal modes {(NNMs)} and their applications in vibration
  theory: An overview.
\newblock {\em Mechanical Systems and Signal Processing}, 11(1):3--22, 1997.

\bibitem{Sarrouy11}
E.~Sarrouy, A.~Grolet, and F.~Thouverez.
\newblock Global and bifurcation analysis of a structure with cyclic symmetry.
\newblock {\em International Journal of Non-Linear Mechanics}, 46(5):727--737,
  2011.

\bibitem{Georgiades09b}
F.~Georgiades and A.~F. Vakakis.
\newblock Passive targeted energy transfers and strong modal interactions in
  the dynamics of a thin plate with strongly nonlinear attachments.
\newblock {\em International Journal of Solids and Structures},
  46(11):2330--2353, 2009.

\bibitem{Vaurigaud11}
B.~Vaurigaud, A.~Ture~Savadkoohi, and C.-H. Lamarque.
\newblock Targeted energy transfer with parallel nonlinear energy sinks. part
  i: Design theory and numerical results.
\newblock {\em Nonlinear Dynamics}, 66(4):763--780, 2011.

\bibitem{Kuether15}
R.~J. Kuether, L.~Renson, T.~Detroux, C.~Grappasonni, G.~Kerschen, and M.~S.
  Allen.
\newblock Nonlinear normal modes, modal interactions and isolated resonance
  curves.
\newblock {\em Journal of Sound and Vibration}, 351:299--310, 2015.

\bibitem{Hill16b}
T.~L. Hill, S.~A. Neild, and A.~Cammarano.
\newblock An analytical approach for detecting isolated periodic solution
  branches in weakly nonlinear structures.
\newblock {\em Journal of Sound and Vibration}, 379(Supplement C):150--165,
  2016.

\bibitem{Peeters11b}
M.~Peeters, G.~Kerschen, and J.~C. Golinval.
\newblock Dynamic testing of nonlinear vibrating structures using nonlinear
  normal modes.
\newblock {\em Journal of Sound and Vibration}, 330(3):486--509, 2011.

\bibitem{Peeters11c}
M.~Peeters, G.~Kerschen, and J.~C. Golinval.
\newblock Modal testing of nonlinear vibrating structures based on nonlinear
  normal modes: Experimental demonstration.
\newblock {\em Mechanical Systems and Signal Processing}, 25(4):1227--1247,
  2011.

\bibitem{Renson16}
L.~Renson, A.~Gonzalez-Buelga, D.~A.~W. Barton, and S.~A. Neild.
\newblock Robust identification of backbone curves using control-based
  continuation.
\newblock {\em Journal of Sound and Vibration}, 367:145--158, 2016.

\bibitem{Peter17}
S.~Peter and R.~I. Leine.
\newblock Excitation power quantities in phase resonance testing of nonlinear
  systems with phase-locked-loop excitation.
\newblock {\em Mechanical Systems and Signal Processing}, 96:139--158, 2017.

\bibitem{Hill16}
T.~L. Hill, P.~L. Green, A.~Cammarano, and S.~A. Neild.
\newblock Fast {B}ayesian identification of a class of elastic weakly nonlinear
  systems using backbone curves.
\newblock {\em Journal of Sound and Vibration}, 360:156--170, 2016.

\bibitem{Lacarbonara16}
W.~Lacarbonara, B.~Carboni, and G.~Quaranta.
\newblock Nonlinear normal modes for damage detection.
\newblock {\em Meccanica}, 51(11):2629--2645, 2016.

\bibitem{Peter15}
S.~Peter, A.~Grundler, P.~Reuss, L.~Gaul, and R.~I. Leine.
\newblock Towards finite element model updating based on nonlinear normal
  modes.
\newblock In {\em Proceedings of the International Modal Analysis Conference
  (IMAC)}, Orlando, USA, 2015.

\bibitem{Song18}
M.~Song, L.~Renson, J.P. No\"el, B~Moaveni, and G.~Kerschen.
\newblock Bayesian model updating of nonlinear systems using nonlinear normal
  modes.
\newblock {\em Structural Control and Health Monitoring}, in review.

\bibitem{Ehrhardt17}
D.~A. Ehrhardt, M.~S. Allen, T.~J. Beberniss, and S.~A. Neild.
\newblock Finite element model calibration of a nonlinear perforated plate.
\newblock {\em Journal of Sound and Vibration}, 392:280--294, 2017.

\bibitem{Ehrhardt16}
D.~A. Ehrhardt and M.~S. Allen.
\newblock Measurement of nonlinear normal modes using multi-harmonic stepped
  force appropriation and free decay.
\newblock {\em Mechanical Systems and Signal Processing}, 76-77:612--633, 2016.

\bibitem{Kerschen13}
G.~Kerschen, M.~Peeters, J.~C. Golinval, and C.~St\'ephan.
\newblock Nonlinear modal analysis of a full-scale aircraft.
\newblock {\em Journal of Aircraft}, 50(5):1409--1419, 2013.

\bibitem{Touze02}
C.~Touz\'e, O.~Thomas, and A.~Chaigne.
\newblock Asymmetric non-linear forced vibrations of free-edge circular plates.
  {P}art 1: Theory.
\newblock {\em Journal of Sound and Vibration}, 258(4):649--676, 2002.

\bibitem{Ehrhardt17sub}
D.A. Ehrhardt, T.~Hill, S.A. Neild, and J.E. Cooper.
\newblock Veering and nonlinear interactions of a clamped beam in bending and
  torsion.
\newblock {\em Journal of Sound and Vibration, in press}.

\bibitem{WaggNeildBook}
D.~Wagg and S.~A. Neild.
\newblock {\em Nonlinear Vibration with Control --- For Flexible and Adaptive
  Structures}.
\newblock Springer, The Netherlands, 2010.

\bibitem{Silva94}
M.~R.~M. Crespo Da~Silva and C.~L. Zaretzky.
\newblock Nonlinear flexural-flexural-torsional interactions in beams including
  the effect of torsional dynamics. i: Primary resonance.
\newblock {\em Nonlinear Dynamics}, 5(1):3--23, 1994.

\bibitem{Westra12}
H.~J.~R. Westra, H.~S.~J. van~der Zant, and W.~J. Venstra.
\newblock Modal interactions of flexural and torsional vibrations in a
  microcantilever.
\newblock {\em Ultramicroscopy}, 120(Supplement C):41--47, 2012.

\bibitem{Polunin17}
P.~M. Polunin and S.~W. Shaw.
\newblock Self-induced parametric amplification in ring resonating gyroscopes.
\newblock {\em International Journal of Non-Linear Mechanics}, 94(Supplement
  C):300--308, 2017.

\bibitem{Kuether15b}
R.~J. Kuether, B.~J. Deaner, J.~J. Hollkamp, and M.~S. Allen.
\newblock Evaluation of geometrically nonlinear reduced-order models with
  nonlinear normal modes.
\newblock {\em AIAA Journal}, 53(11):3273--3285, 2015.

\bibitem{RensonISMA2016}
L.~Renson, D.A. Ehrhardt, D.A.W. Barton, S.A. Neild, and J.E. Cooper.
\newblock Connecting nonlinear normal modes to the forced response of a
  geometric nonlinear structure with closely spaced modes.
\newblock In {\em Proceedings of the ISMA}, Leuven, BE, 2016.

\bibitem{RensonReview}
L.~Renson, G.~Kerschen, and B.~Cochelin.
\newblock Numerical computation of nonlinear normal modes in mechanical
  engineering.
\newblock {\em Journal of Sound and Vibration}, 364:177--206, 2016.

\bibitem{Lacarbonara05}
W.~Lacarbonara, H.~N. Arafat, and A.~H. Nayfeh.
\newblock Non-linear interactions in imperfect beams at veering.
\newblock {\em International Journal of Non-Linear Mechanics}, 40(7):987--1003,
  2005.

\bibitem{Wright99}
J.~R. Wright, J.~E. Cooper, and M.~J. Desforges.
\newblock Normal mode force appropriation - theory and application.
\newblock {\em Mechanical Systems and Signal Processing}, page~24, 1999.

\bibitem{Zapico13}
J.~L. Zapico-Valle, M.~Garcia-Di\'eguez, and R.~Alonso-Camblor.
\newblock Nonlinear modal identification of a steel frame.
\newblock {\em Engineering Structures}, 56(0):246--259, 2013.

\bibitem{Peeterse11c}
M.~Peeters, G.~Kerschen, and J.~C. Golinval.
\newblock Modal testing of nonlinear vibrating structures based on nonlinear
  normal modes: Experimental demonstration.
\newblock {\em Mechanical Systems and Signal Processing}, 25(4):1227--1247,
  2011.

\bibitem{Hill2014b}
T~L Hill, A~Cammarano, S~A Neild, and D~J Wagg.
\newblock An analytical method for the optimisation of weakly nonlinear
  systems.
\newblock {\em Proceedings of EURODYN 2014}, pages 1981--1988, 2014.

\bibitem{Hill15}
T.~L. Hill, A.~Cammarano, S.~A. Neild, and D.~J. Wagg.
\newblock Interpreting the forced response of a two-degree-of-freedom nonlinear
  oscillator using backbone curves.
\newblock {\em Journal of Sound and Vibration}, 349:276--288, 2015.

\bibitem{TouzeCISM}
C.~Touz\'e.
\newblock {\em Normal form theory and nonlinear normal modes: Theoretical
  settings and applications}, pages 75--160.
\newblock Springer Vienna, Vienna, 2014.

\end{thebibliography}

\end{document}